\documentclass[12pt]{article}
\usepackage{amsmath}
\usepackage{amsbsy}
\usepackage{amssymb}
\usepackage{latexsym}
\usepackage{mathrsfs}
\usepackage[usenames]{color}
\usepackage{fullpage}
\usepackage[colorlinks=true,citecolor=OliveGreen,linkcolor=BrickRed,urlcolor=BrickRed]{hyperref}
\usepackage{breakurl}

\def\bbC{{\mathbb C}}

\def\bbZ{{\mathbb Z}}
\def\bbN{{\mathbb N}}

\def\bbR{{\mathbb R}}
\def\bbQ{{\mathbb Q}}

\def\bfa{{\boldsymbol A}}
\def\bfb{{\boldsymbol B}}
\def\bfc{{\boldsymbol C}}
\def\bfd{{\boldsymbol D}}
\def\bfe{{\boldsymbol E}}

\def\bfk{{\boldsymbol K}}

\def\bfu{{\boldsymbol U}}
\def\bfx{{\boldsymbol X}}
\def\bfv{{\boldsymbol V}}
\def\bfw{{\boldsymbol W}}

\def\bfone{{\boldsymbol 1}}

\def\a{{\cal A}}

\def\h{{\cal H}}
\def\k{{\cal K}}

\def\n{{\cal N}}

\def\calu{{\cal U}}

\newtheorem{thm}{Theorem}[section]
\newtheorem{lem}[thm]{Lemma}

\newtheorem{cor}[thm]{Corollary}
\newtheorem{pro}[thm]{Proposition}

\newtheorem{question}[thm]{Question}

\def\verts{{V}}
\def\ord{{<}}

\def\Nind{N_{\rm ind}}
\def\Nord{N_{\rm ord}}
\def\Nindr{N_{\rm ind}^{\rm res}}
\def\Nordr{N_{\rm ord}^{\rm res}}
\def\GG{{\bf G}}
\def\HH{{\bf H}}
\def\YY{{\bf Y}}
\def\XX{{\bf X}}
\def\forb{{\h}}
\def\Hilb{{\mathscr H}}
\def\E{{\bf E}}
\def\P{{\bf P}}
\def\bfPhi{{\boldsymbol\Phi}}
\def\bford{{\boldsymbol\ord}}
\def\Seq#1{\langle #1 \rangle}
\def\ceil#1{\lceil #1 \rceil}
\def\dfnterm#1{\textit{\textbf{#1}}}
\newcommand{\TV}{\mathrm{TV}}
\def\Aut{{\rm Aut}}
\def\st{\,;\;}
\def\arXiv#1{\url{http://www.arxiv.org/abs/#1}}
\def\restrict{\mathord{\upharpoonright}}   
\def\gp{{\Gamma}}
\def\fra{Fra\"{i}ss\'e}
\def\symdiff{\mathbin{\triangle}}  
\def\II#1{\bfone_{#1}}

\def\thmenv#1#2#3{\begin{#1} \label{#1:#2} #3 \end{#1}}

\def\procl#1.#2 #3\endprocl{%
       \ifx#1t\thmenv{thm}{#2}{#3}\fi
       \ifx#1l\thmenv{lem}{#2}{#3}\fi
       \ifx#1p\thmenv{pro}{#2}{#3}\fi
       \ifx#1c\thmenv{cor}{#2}{#3}\fi
       \ifx#1d\thmenv{defn}{#2}{#3}\fi
       \ifx#1g\thmenv{conj}{#2}{#3}\fi
       \ifx#1q\thmenv{question}{#2}{#3}\fi
       \ifx#1r\thmenv{remark}{#2}{{\rm #3}}\fi
    }%

\def\rref#1.#2/{%
      \ifx #1sSection~\ref{s.#2}\fi
      \ifx #1tTheorem~\ref{thm:#2}\fi  
      \ifx #1lLemma~\ref{lem:#2}\fi 
      \ifx #1cCorollary~\ref{cor:#2}\fi 
      \ifx #1pProposition~\ref{pro:#2}\fi 
      \ifx #1dDefinition~\ref{defn:#2}\fi
      \ifx #1gConjecture~\ref{conj:#2}\fi 
      \ifx #1qQuestion~\ref{question:#2}\fi 
      \ifx #1rRemark~\ref{remark:#2}\fi 
      \ifx #1aAppendix~\ref{a.#2}\fi 
      \ifx #1fFigure~\ref{f.#2}\fi
      \ifx #1e(\ref{e.#2})\fi
      \ifx #1b\cite{#2}\fi
        }

\def\rlabel #1 #2{\begin{equation} \label{#1} #2 \end{equation}}

\def\proof{\medbreak\noindent{\it Proof.\enspace}}

\def\proofof #1.#2 {\medbreak\noindent
     {\it Proof of \rref #1.#2/.}\enspace}

\def\Qed{\hfill$\dashv$\medbreak}

\def\bsection#1#2{\bigbreak\section{#1}\label{#2}}

\begin{document}
\title{\bf Random Orderings and Unique Ergodicity of Automorphism Groups}
\author{Omer Angel\footnote{Research partially supported by NSERC and the
Sloan Foundation}, Alexander S. Kechris\footnote{Research partially
supported by NSF grant DMS-0968710.}, and Russell Lyons\footnote{Research
partially supported by NSF grant DMS-1007244 and Microsoft Research.}}
\date{9 August 2012}
\maketitle

\begin{abstract}
We show that the only random orderings of finite graphs
that are invariant under isomorphism and induced
subgraph are the uniform random orderings. We show how this implies the
unique ergodicity of the automorphism group of the random graph.
We give similar theorems for other structures, including, for example, metric spaces.
These give the first examples of uniquely ergodic groups, other than
compact groups and extremely amenable groups, after Glasner
and Weiss's example of the group of all permutations of the integers. We also contrast these results to those for certain special classes of graphs and metric spaces in which such random orderings can be found that are not uniform.
\end{abstract}

\section{Introduction}
\subsection{Random Orderings}

Consider the class of finite graphs, by which we mean simple graphs, i.e.,
without loops or multiple edges.
Is there any way to distinguish among the vertices of a finite graph
in a way that is preserved by isomorphism and by taking induced subgraphs?
To make this question more precise, consider random linear (total) 
orderings of vertices of finite graphs.
That is, for each graph $G = \langle V, E\rangle$, let $\mu_G$ be a probability measure
on the $|V|!$ linear orderings of $V$.
Suppose that the collection of measures $\mu_G$ is
\dfnterm{consistent}, meaning that it satisfies two properties:
\begin{enumerate}
\item[i)]  If $\phi : G \to G'$ is a graph isomorphism, then $\phi_* \mu_G
= \mu_{G'}$. Here, $\phi_*$ denotes the push-forward map induced by $\phi$;
more precisely, $\phi_*$ is the push-forward of the map $(\phi, \phi)$ on
orders.
\item[ii)] If $H$ is an induced subgraph of $G$, then $\mu_G$ induces
$\mu_H$ by restriction. In other words, if $\phi$ is the restriction map of
orderings of $V(G)$ to orderings of $V(H)$, then $\phi_* \mu_G = \mu_H$.
\end{enumerate}
We shall refer to the family $(\mu_G)_G$ as a \dfnterm{consistent random ordering} (for the class of finite graphs). Note that property (i) by itself guarantees that for a complete graph $G$, as well as
an empty graph, $\mu_G$ must be the uniform measure. If $(\mu_G)_G$ is
consistent, must $\mu_G$ be the
uniform measure on all linear orderings of $V(G)$ for all $G$, or is there a more interesting consistent assignment of random orderings?

For example, if instead of the collection of all finite graphs, we
considered consistent random orderings only of finite graphs that are
paths, then there is clearly another choice: orient the path in one of the
two ways at random, with probability 1/2 each, and use the orientation to
define the naturally associated ordering.

It is much harder to define a non-uniform consistent random ordering on
the class of finite graphs all of whose components are paths, but it can be
done.
What about graphs with a given bound on their degrees, or other
classes of graphs? What about finite hypergraphs or finite metric spaces?
In fact, such questions can be asked in great generality for classes of finite
structures in a given language in the sense of model theory.

We show in this paper the following.

\procl t.ex
  The only consistent random ordering for the class of finite graphs is the
  uniform ordering. The same holds for the classes of $K_n$-free graphs,
  $r$-uniform hypergraphs, and metric spaces with (non-zero) distances in a given
  additive subsemigroup of\/ $\bbR^+$. 
\endprocl

This is proved in Sections \ref{s.quant}--\ref{s.eucl}, where many other
such examples are given, including classes of hypergraphs with forbidden
configurations. In these sections, we also discuss several examples of
classes of metric spaces and graphs for which the opposite happens, i.e.,
there are non-uniform consistent random orderings, including the class of
Euclidean metric spaces and the class of bounded degree graphs, which can
be proved by using a random projection method suggested to us by Leonard
Schulman.

Furthermore, we obtain quantitative versions of \rref t.ex/.  Let
$d_\TV(\mu,\nu) := \frac12\|\mu-\nu\|_1 = \max_A |\mu(A) - \nu(A)|$ denote
the total variation distance between probability measures.  We show the
following:

\procl t.approx1
Let $2 \le r \le k \le n$ be integers.
There is a constant $C = C(k, r)$ with the following property.
  Let $(\mu_G)_G$ be a consistent random ordering defined on all $r$-uniform
  hypergraphs $G$ of
  size at most $n$, and let $H$ be a $r$-uniform
  hypergraph of size $k$. If $\nu_H$ is the uniform ordering on $H$, then
  $d_\TV(\mu_H,\nu_H) \le C \sqrt{\log n \over n^{r-1}}$.
\endprocl

In the case of graphs, we construct in Section~\ref{s.lower} a random
ordering on graphs of size at most $n$ (or even with degrees bounded by
$n$) such that for some $H$, $d_\TV(\mu_H,\nu_H) \geq C/n$. A similar
construction appears to give $d_\TV(\mu_H,\nu_H) \geq C/n^{r-1}$ for
$r$-uniform hypergraphs.

\begin{question}
  What is the largest possible total variation distance from the uniform
  ordering of a random ordering on graphs (or hypergraphs) of size at most $n$?
\end{question}

\subsection{Unique Ergodicity}

The reason for our attention to these questions, beyond their intrinsic
interest, is that, in certain circumstances, they provide a way to prove unique ergodicity results for groups.
In order to explain this context, we need to review some concepts and results concerning the model theory and combinatorics of classes of finite structures.

The general setting for our ergodicity results is the following.
A (first-order) \dfnterm{language} $L$ consists of families $ (R_i)_{i\in
I}, (f_j)_{j\in J}$ of relation symbols $R_i$ and function symbols $f_j$
with associated arities $m_i \geq 1$ and $n_j\geq 0$. A \dfnterm{structure}
for this language, or \dfnterm{$L$-structure}, is an object of the form $$\bfa
=  \langle A, (R_i^\bfa)_{i\in I}, (f_j^\bfa)_{j\in J}\rangle,$$where $A$
is a non-empty set, called the \dfnterm{universe} of the structure,
$R_i^\bfa\subseteq A^{m_i}$ and $f_j^\bfa\colon A^{n_j}\to A$, where if
$n_j = 0$, it is understood that $f_j^\bfa$ is just an element of $A$. The
cardinality of the structure is the cardinality of its universe $A$.
\dfnterm{In this paper, all languages and structures will be countable.} For brevity, and when there is no danger of confusion, we sometimes omit the superscripts.

As an example, if $L = \{E\}$ with $E$ a binary relation symbol, then the
class of graphs is the class of $L$-structures $\bfa$ for which $E^\bfa$ is
symmetric and irreflexive. A metric space $\langle X, d\rangle$ can also be
viewed as a structure $\bfx = \langle X, (R_q^\bfx)_{q\in\bbQ^+}\rangle$ in
the language with binary relation symbols $(R_q)_{q\in\bbQ^+}$, where
$R_q^\bfx (x,y) \iff d(x,y) < q$.

Although this will be our standard notation when we discuss abstract
structures, we shall keep the more traditional notation (mainly in font type) for specific structures, like graphs, hypergraphs, metric spaces, etc., that we used earlier and that we also use in Sections \ref{s.quant}--\ref{s.eucl}.

A class $\k$ of finite $L$-structures is called a \dfnterm{Fra\"{i}ss\'e class} if it contains 
structures of arbitrarily large (finite) cardinality, is countable
(in the sense that it contains only countably many isomorphism types)
and satisfies the following:
\begin{enumerate}
\item[i)] \dfnterm{Hereditary Property}: If $\bfb\in\k$ and $\bfa$ can
be embedded in $\bfb$, then $\bfa\in\k$.
\item[ii)] \dfnterm{Joint Embedding Property}: If $\bfa ,\bfb\in\k$,
there is $\bfc\in\k$ such that $\bfa ,\bfb$ can be embedded in $\bfc$.
\item[iii)] \dfnterm{Amalgamation Property}: If $\bfa ,\bfb ,\bfc
\in\k$ and $f\colon\bfa\to\bfb ,g\colon\bfa\to\bfc$ are embeddings,
there is $\bfd\in\k$ and embeddings $r\colon\bfb\to\bfd$, $s\colon\bfc
\to\bfd$ such that $r\circ f=s\circ g$.
\end{enumerate}

Throughout this paper embeddings and substructures will be understood in
the usual model theoretic sense (see, e.g., Hodges \cite{Ho}, page 5);
e.g., for graphs embeddings are induced embeddings, i.e., isomorphisms onto
induced subgraphs. To be precise, given a language $L$ consisting of
families $ (R_i)_{i\in I}$ and $(f_j)_{j\in J}$ of relation symbols $R_i$ and
function symbols $f_j$ with associated arities $m_i \geq 1$ and $n_j\geq
0$, an $\dfnterm{embedding}$ of an $L$-structure $\bfa$ into an
$L$-structure $\bfb$ is an injection $\phi\colon A\to B$ such that for any
$R_i$ and $a_1, \dots , a_{m_i} \in A$, we have $R_i^\bfa (a_1, \dots ,
a_{m_i})\iff R_i^\bfb \big(\phi(a_1), \dots , \phi(a_{m_i})\big)$ and for any $f_j$
and any $a_1, \dots , a_{n_j} \in A$, we have $\phi\big( f_j^\bfa (a_1, \dots ,
a_{n_j})\big) =f_j^\bfb \big(\phi(a_1), \dots , \phi(a_{n_j})\big)$. If the identity is
such an embedding, then we say that $\bfa$ is a \dfnterm{substructure} of
$\bfb$.

We recall the following results of \cite{Fr} (see also \cite[Section 7.1]{Ho}):

If $\k$ is a Fra\"{i}ss\'e class, there is a unique, up to isomorphism,
countably infinite structure $\bfk$ that is \dfnterm{locally finite} (i.e.,
each finite subset of $K$ is contained in a finite substructure of $\bfk$),
\dfnterm{ultrahomogeneous} (i.e., any isomorphism between finite
substructures can be extended to an automorphism of the structure) and is
such that, up to isomorphism, its finite substructures are exactly those in
$\k$.
We call this the \dfnterm{Fra\"{i}ss\'e limit} of $\k$, in symbols
\[
\bfk = {\rm Flim}(\k ).
\]

Conversely, if $\bfk$ is a countably infinite structure that is locally finite and ultrahomogeneous,
then its \dfnterm{age}, Age$(\bfk )$, i.e., the class $\k$ of all its
finite substructures, up to isomorphism, is a Fra\"{i}ss\'e class.
Such structures $\bfk$ are called \dfnterm{Fra\"{i}ss\'e structures}.
Thus there is a canonical bijection $\k\mapsto{\rm Flim}(\k )$,
$\bfk\mapsto{\rm Age}(\bfk )$ between Fra\"{i}ss\'e classes and structures.

We are interested in the set of invariant Borel probability measures for continuous actions of the automorphism group
$\gp:={\rm Aut}(\bfk )$, viewed as a topological group under the
pointwise convergence topology.  We note that the groups Aut$(\bfk )$,
for \fra\ structures $\bfk$ as above, are (up to topological group isomorphism) exactly the closed subgroups of the infinite
symmetric group $S_\infty$, i.e, the group of permutations of $\bbN$, again
with the pointwise convergence topology (see, e.g., \cite{KPT}). Up to
topological group isomorphism they are also the same as the non-archimedean
Polish groups, where a topological group is called
\dfnterm{non-archimedean} if it admits a basis at the identity consisting
of open subgroups (see \cite[1.5.1]{BK}).

Consider now amenability properties of such groups. At one end of the spectrum, there are many examples of automorphism groups
that satisfy a very strong form of amenability, i.e., they are
\dfnterm{extremely amenable}.  This means that every continuous action of
such a group on a (non-empty) compact Hausdorff space, i.e., a
\dfnterm{$\gp$-flow}, has a fixed point.  As was shown in \cite[4.7]{KPT},
$\gp={\rm Aut}(\bfk )$ has this property exactly when  $\bfk$ is an
\dfnterm{order} Fra\"{i}ss\'e structure, i.e., a Fra\"{i}ss\'e structure in
which one of the relations is a linear ordering, such that Age$(\bfk )$
satisfies the Ramsey Property (RP); see \cite[Section 3]{KPT}. We also
discuss the RP in  \rref s.S2/ below.  Extensive lists of extremely
amenable automorphism groups are discussed in \cite[Section 6]{KPT}.

Next there are automorphism groups that are \dfnterm{amenable} (i.e.,
every $\gp$-flow has an invariant Borel probability measure) but not
extremely amenable.  The most extensive list, in our framework,
of such examples 
arises in the context of the Hrushovski property.  Given a class
$\k$ of finite structures in a given language, we say that $\k$
is a  \dfnterm{Hrushovski class} if for any $\bfa\in\k$
and any (partial) isomorphisms $\varphi_i\colon \bfb_i\to \bfc_i$, 
$1\leq i\leq k$, where $\bfb_i$, $\bfc_i$ are substructures of $\bfa$,
there is $\bfb\in\k$ containing $\bfa$ as a substructure such that
each $\varphi_i$ can be extended to an automorphism $\psi_i$ of $\bfb$,
$1\leq i\leq k$.  It is shown in \cite[6.4]{KR} that if $\k$ is a Fra\"{i}ss\'e
class with $\bfk ={\rm Flim}(\k )$, then $\k$ is a Hrushovski class iff $\gp=
{\rm Aut}(\bfk )$ is \dfnterm{compactly approximable}, i.e.,  there is an
increasing sequence $\Delta_0\subseteq \Delta_1\subseteq\cdots$ of compact
subgroups of $\gp$ with $\overline{\bigcup_n\Delta_n}=\gp$.  Calling the
Fra\"{i}ss\'e limit $\bfk$ of a Hrushovski class $\k$ a \dfnterm{Hrushovski structure}, we thus see
that the automorphism groups of Hrushovski structures are compactly
approximable, whence amenable.  Examples of Hrushovski classes of finite structures include the
following: pure sets, graphs (Hrushovski \cite{Hr}), 
$r$-uniform hypergraphs, $K_n$-free graphs
(Herwig \cite{He1,He2}), metric spaces with distances in a countable additive subsemigroup of $\bbR^+$ (Solecki \cite{So}),
finite-dimensional vector spaces over a fixed finite
field, etc.

Finally, at the other end of the spectrum, there are also
automorphism groups $\gp={\rm Aut}(\bfk )$ that are not amenable.
These include, among others, the automorphism groups of the countable
atomless Boolean algebra, the random poset, and
the random distributive lattice (see \cite{KS}).

A characterization of
the amenability of $\gp={\rm Aut}(\bfk )$ in terms of combinatorial
properties of the Fra\"{i}ss\'e class Age$(\bfk )$ was obtained by
Moore and Tsankov, see \cite[6.1]{Mo}.  Another characterization, in a special
case, is in \rref p.amen/ below.

In this paper we are interested in the ergodic theory of the flows of
automorphism groups, in particular, in the potential phenomenon
of unique ergodicity. \dfnterm{Below, measure means Borel probability
measure}.

Let $\gp$ be a topological group. We say that
a $\gp$-flow is \dfnterm{uniquely ergodic} if it has a unique invariant
measure (which therefore must be ergodic).  We say that the group 
$\gp$ is \dfnterm{uniquely ergodic} if every \dfnterm{minimal} $\gp$-flow
(i.e., one all of whose orbits are dense) is uniquely ergodic.
(The assumption of minimality is clearly necessary, as a given
$\gp$-flow may have in general many minimal subflows that are
pairwise disjoint, and by amenability each will support an
invariant measure.) Clearly every extremely amenable Polish group is uniquely ergodic and so is every compact Polish group.

However, this property is never realized in the realm of infinite countable (discrete) groups, as follows from results of Benjamin Weiss \cite{W}. Weiss also believes that this extends to non-compact, Polish locally compact groups, although this has not been checked in detail yet.

It is important here to review the concept of the universal minimal 
flow of a topological group, $\gp$. It is a classical result in
topological dynamics that every topological group $\gp$ admits a
unique, up to isomorphism of $\gp$-flows, minimal flow $M(\gp)$, called
its \dfnterm{universal minimal flow}, such that all other minimal
$\gp$-flows are factors of it (see, e.g.,
\cite[Section 1]{KPT}). Recall that  a $\Gamma$-flow $Y$ is a
\dfnterm{factor} of a $\Gamma$-flow $X$ if there is a continuous surjection
$\phi\colon X\to Y$ that is a $\Gamma$-map, i.e., $\phi (\gamma \cdot x) =
\gamma \cdot \phi (x)$ for every $x\in X$ and $\gamma\in \Gamma$. If such a map
is a bijection, we call it an \dfnterm{isomorphism}. Since every $\gp$-flow
contains a minimal subflow, the group $\gp$ is amenable iff $M(\gp)$ has an
invariant measure. Also, it can be seen that $\gp$ is uniquely ergodic iff
$M(\gp)$ is uniquely ergodic (see \rref p.Char/).

When $\gp$ is compact, $M(\gp)= \gp$ (with the left translation action) and
when $\gp$ is extremely amenable, $M(\gp)$ is trivial, i.e., a singleton, but in general $M(\gp)$ is a very complicated object that is difficult to
``compute''.  For example, when $\gp$ is infinite countable, $M(\gp)$ is a ``big''
space of ultrafilters on $\gp$, and more generally when $\gp$ is non-compact,  locally
compact, $M(\gp)$ is not metrizable (see \cite[A2.2]{KPT}).  However, over the last fifteen years or so, there have been many examples of explicit descriptions of
non-trivial metrizable universal minimal flows: see Pestov \cite{Pe1},
Glasner and Weiss
\cite{GW1, GW2} and, in the case
of automorphism groups of Fra\"{i}ss\'e structures, \cite{KPT}.

As was shown in \cite{GW1}, where the universal minimal flow of
$S_\infty$ was computed, $S_\infty$ is uniquely ergodic.  This seems to be
the first example found of such a group that is neither extremely amenable
nor compact. We shall show in this paper that many other amenable automorphism groups,
for which we can compute a metrizable universal minimal flow, are also uniquely ergodic, so this appears as a general phenomenon. This will be a consequence of the uniqueness results for consistent random orderings, like \rref t.ex/, that we prove in Sections \ref{s.quant}--\ref{s.metric}.

To see this connection, we need to explain briefly the method by which
universal minimal flows are calculated in \cite{KPT}. Details and precise
definitions are given in \rref s.S2/. Given a Fra\"{i}ss\'e class $\k$ with
Fra\"{i}ss\'e limit $\bfk$, one appropriately assigns to each structure $\bfa \in \k$ a
collection of so-called {admissible} (linear) orderings on the universe $A$
of $\bfa$ to obtain a new class $\k^*$  of expanded, ordered structures of
the form $\langle \bfa, <\rangle$, where $\bfa \in \k$ and $<$ is an
admissible ordering on $\bfa$. In many cases, such as for graphs,
every linear ordering on a given structure is admissible, but in other
cases, one has to take a more restricted collection of linear orderings.
For example, take as $\k$ the Fra\"{i}ss\'e class of bipartite graphs with
distinguished parts, i.e., structures of the form $\langle V, E, A_1,
A_2\rangle$, where $\langle V, E\rangle$ is a graph and $A_1, A_2$ is a
partition of $V$ such that there are no edges within each $A_i$. Then for
each such graph, the relevant collection of admissible orderings consists of all orderings $<$ such that $A_1 < A_2$.

If an appropriate such class $\k^*$ of expanded structures as
above can be found that satisfies various structural conditions, including
Ramsey properties, we shall call $\k^*$ a ``companion" of $\k$. In that
case, it is shown in \cite{KPT} that the universal minimal flow of
$\gp={\rm Aut}(\bfk )$ is obtained as follows. Denote by $X_{\k^*}$ the
space of all orderings $<$ on the universe $K$ of $\bfk$ that are admissible (relative to
$\k^*$) in the sense that for every finite substructure $\bfa$ of $\bfk$,
the restriction of $<$ to $A$ is admissible (for $\bfa$). This is a
compact metrizable space on which $\gp$ acts continuously in the obvious
way, and it turns out to be the universal minimal flow for $\gp$. From this one can see that the existence of an invariant measure on $X_{\k^*}$ (i.e., the amenability of $\gp$) is equivalent to the existence of a consistent random admissible ordering and unique ergodicity of $\gp$ is equivalent to the uniqueness of a consistent random admissible ordering. Using this and the results in Sections \ref{s.quant}--\ref{s.metric} concerning uniqueness of consistent random orderings, we obtain, in Sections \ref{s.ot}--\ref{s.qo}, many new examples of uniquely ergodic automorphism groups. A sample is included in the following theorem, which we state after we introduce some terminology.

If $L = (R_i)_{i\in I}$ is a finite relational language with $R_i$ of arity
$m_i\geq 2$, then a \dfnterm{hypergraph of type $L$} is an $L$-structure
$\bfa$ in which each $R_i^\bfa$ gives an $m_i$-uniform hypergraph. A
hypergraph of type $L$ is  called \dfnterm{irreducible} if it has at least
two vertices and every two-element subset of the vertices belongs to some
hyperedge. Given a class $\a$ of irreducible hypergraphs of type $L$, a
hypergraph of type $L$ is \dfnterm{$\a$-free} if it contains no (induced) copy of a
structure in $\a$. If we choose $L$ to have only one relation symbol of
arity $r$ and $\a=\emptyset$, then we obtain the class of $r$-uniform
hypergraphs (graphs if $r=2$), and if we choose $\a=\emptyset$, we obtain
the class of all hypergraphs of type $L$. If we choose $L$ to have only one
symbol of arity $2$ and $\bfa = \{ K_n\}$, then we obtain the class of
$K_n$-free graphs, where $K_n$ is the complete graph on $n$ vertices.

The \dfnterm{random $\a$-free hypergraph} of a given type $L$ is the
Fra\"{i}ss\'e limit of the class of $\a$-free hypergraphs of type $L$ (thus
by choosing $L$ and $\a$ appropriately, this includes the case of the random
graph, random $K_n$-free graph, random $r$-uniform hypergraph, etc.). The
\dfnterm{Urysohn space} $\bfu_S$, where $S$ is a countable additive
subsemigroup of $\bbR^+$, is the Fra\"{i}ss\'e limit of the class of finite
metric spaces with distances in $S$. 

\procl t.ex1
The automorphism groups of the following structures are uniquely ergodic: the (countably) infinite-dimensional vector space over a given finite field, the random $\a$-free hypergraph of a given type,  and  $\bfu_S$.
\endprocl

Restricting ourselves to automorphism groups of Hrushovski
structures, which provide the most prominent examples of amenable (but
not extremely amenable) groups, we shall find that unique ergodicity of
Aut$(\bfk )$, with $\bfk$ a Hrushovski structure, is equivalent to a
combinatorial property of Age$(\bfk )$, very much in the spirit of
\cite{KPT}.  In fact, rather interestingly, if $\k := \textrm{Age}(\bfk)$ admits a companion $\k^*$ as above, then it turns out that unique
ergodicity is exactly equivalent to a quantitative version of what is
called the
``ordering property''. The ordering property is
a key ingredient of the Ramsey theory of
classes of finite structures that is instrumental in the computation 
of universal minimal flows in \cite[7.5]{KPT}.  We discuss this in
\rref s.hru1/ below.

In \rref s.sup/, we show that for certain automorphism groups, including
those of the random $\a$-free uniform hypergraph of a given type and
of $\bfu_S$, every minimal action not only has a unique invariant measure, but
also this measure concentrates on a single comeager orbit. This was earlier
proved for the group $S_{\infty}$ by Glasner and Weiss \cite{GW1}. 

Finally in the last \rref s.problems/, we discuss some open problems arising from the work in this paper.

\bsection{Graphs and Uniform Hypergraphs}{s.quant}

In this section, we prove that the only consistent random ordering on the
class of all finite graphs is the uniform ordering. In fact, we prove the
same for hypergraphs.  Recall that an \dfnterm{$r$-uniform hypergraph} is a
pair $G = \langle V, E\rangle$, where $E \subseteq {V \choose r}$ is a
collection of subsets of $V$ of cardinality $r$; the elements of $E$ are
called \dfnterm{hyperedges}.  The case $r = 2$ is the case of graphs.  The
\dfnterm{size} of $G$ is defined to be the cardinality of $V$.  If $G =
\langle V, E\rangle$ is a hypergraph and $V' \subseteq V$, then the
hypergraph \dfnterm{induced} on $V'$ by $G$ equals $\big\langle V', E \cap
{V' \choose r}\big\rangle$. Note that hyperedges intersecting $V'$ that are not
contained in $V'$ are discarded.




The way we prove \rref t.approx1/ is via the following general
principle.
Let $\Nind(H, G)$ denote the number of embeddings of $H$ in $G$, i.e.,
the number of isomorphisms $\phi : H \to H'$ such that $H'$ is an induced
hypergraph in $G$. (Up to symmetries, this is the number of induced
subgraphs of $G$ that are isomorphic to $H$.)
Given a pair of orderings $\ord_G$ of $V(G)$ and $\ord_H$ of $V(H)$, let
$\Nord(H, G)$ denote the number of ordered embeddings of $H$ in
$G$, i.e.,
the number of embeddings $\phi : H \to G$ such that 
$\phi^*(\ord_G) = \ord_H$.
Here, $\phi^*$ denotes the pull-back map induced by $\phi$, i.e., $x<_{\phi^*(\ord_G)} y \iff \phi (x) <_G \phi (y).$

\procl l.TVbd
Let $k \ge r \ge 2$ be integers.
Let $G$ be an $r$-uniform hypergraph 
and $H$ be an $r$-uniform hypergraph on $k$ vertices such that $\Nind(H, G)
> 0$.
Suppose $\delta$ is such that
for every pair of orderings $\ord_G$ of $V(G)$ and $\ord_H$ of $V(H)$, 
\rlabel e.eachbd
{\left| {\Nord(H, G) \over \Nind(H, G)} - {1 \over k!} \right|
\le \delta
\,.
}
Let $\mu_G$ and $\mu_H$ be random orderings on $G$ and $H$, respectively.
Suppose that every embedding $\phi$ of $H$ in $G$ satisfies $\phi^* \mu_G =
\mu_H$.
Then $d_\TV(\mu_H,\nu_H) \le \delta k!/2$, where $\nu_H$ is the uniform ordering on $H$.
\endprocl

\proof
Fix $\ord_H$.
Choose $\bford_G$ at random according to $\mu_G$ and choose an embedding
$\bfPhi$ of $H$ in $G$ uniformly at random.
Let $A$ be the event that the restriction of $\bford_G$ to the image of
$\bfPhi$ equals $\bfPhi_*(\ord_H)$.
Since $\phi^* \mu_G = \mu_H$ for every $\phi$, we have
$\P(A \mid \bfPhi = \phi) = \mu_H(\ord_H)$,
whence averaging over $\phi$ gives $\P(A) = \mu_H(\ord_H)$.
We can rewrite the assumption \rref e.eachbd/
as $|\P(A \mid \bford_G = \ord_G) - 1/k!| \le \delta$ for each $\ord_G$;
averaging over $\ord_G$ gives $|\P(A) - 1/k!| \le \delta$. 
That is, $|\mu_H(\ord_H) - 1/k!| \le \delta$.
Finally, summing over all orderings $\ord_H$ gives the bound.
\Qed

Clearly \rref t.approx1/ follows from \rref l.TVbd/ and the following
result. Write $I(n,k,r) := (n)_k 2^{-{k \choose r}}$, where
$(n)_k := n(n-1) \dots (n-k+1)$ is the number of 1-1 maps from $\{1,\dots,
k\}$ to $\{1,\dots, n\}$.

\procl t.quant
Let $2 \le r \le k \le n$ be integers.
There is a constant $C = C(k, r)$ with the following property.
For every $r$-uniform hypergraph $H$ on $k$ vertices,
there exists an $r$-uniform hypergraph $G$ on $n$ vertices such that
$$
\left|{\Nind(H, G) \over I(n, k, r)} - 1 \right|
<
C \sqrt{\log n \over n^{r-1}}
$$
and for every pair of orderings $\ord_G$ of $V(G)$ and $\ord_H$ of $V(H)$, 
\rlabel e.want
{\left| {\Nord(H, G) \over \Nind(H, G)} - {1 \over k!} \right|
<
C \sqrt{\log n \over n^{r-1}}
\,.
}
\endprocl

The proof of \rref t.quant/ uses the following classical inequality of
McDiarmid \rref b.McD/, known as the \dfnterm{bounded-differences inequality}:

\procl t.McD
Let ${\bf Z} := \Seq{Z_1, \ldots, Z_n}$, where $Z_i$ are
independent random variables, and
$f(z_1, \ldots, z_n)$ be a real-valued function such that 
$$
|f({\bf z}) - f({\bf z'})| \le a_i
$$
when the vectors ${\bf z}$ and ${\bf z'}$ differ only in the $i$th
coordinate.
Write $\zeta := \E\big[f({\bf Z})\big]$.
Then for all $L > 0$,
$$
\P \big [ | f({\bf Z}) - \zeta | \geq L \big ] 
\leq 2\exp \left(- {2 L^2 \over \sum_{i=1}^n a_i^2} \right) \, .
$$
\endprocl

\proofof t.quant
Let $\GG$ be a uniformly random $r$-uniform hypergraph on $n$ fixed
vertices, $V$ (so that each hyperedge is present with probability
$1/2$). Note that
\[
I(n,k,r) = \E\big[\Nind(H, \GG)\big] \,.
\]
Define
\[
f(\GG) := \frac{\Nind(H,\GG)}{I(n,k,r)} \,,
\]
which we consider as a function of the $\binom{n}{r}$ variables indicating
the presence of each possible hyperedge. The addition or removal a single
hyperedge to $\GG$ changes $\Nind(H,\GG)$ by at most $(k)_r (n-r)_{k-r}$,
and so $f$ satisfies the conditions of \rref t.McD/ with $a_i=c_1 n^{-r}$,
where we shall denote by $c_j$ intermediate constants that
depend on $k$ and $r$, but not on $n$.  It follows that
\[
\P\left[ \left| \frac{\Nind(H,\GG)}{I(n,k,r)} - 1 \right | \geq D \right]
\leq
2\exp\left\{ \frac{-2D^2}{\binom{n}{r} (c_1 n^{-r})^2} \right\}
\le
2e^{-c_2 D^2 n^r}
\,.
\]

Similarly, for any fixed orderings $\ord_H$ and $\ord_V$, we have
$\E\big[\Nord(H,\GG)] = I(n,k,r)/k!$.
We apply \rref t.McD/ to 
\[
\GG \mapsto \frac{\Nord(H,\GG)}{I(n,k,r)} \,.
\]
Here, adding or removing a single hyperedge changes $\Nord(H,\GG)$ by at
most $\binom{k}{r} (n-r)_{k-r}$, so as above,
\[
\P\left[ \left| \frac{\Nord(H,\GG)}{I(n,k,r)} - \frac{1}{k!} \right | \geq
  D \right] 
\leq
2\exp\left\{ \frac{-2D^2}{\binom{n}{r} (c_3 n^{-r})^2} \right\}
\le
2e^{-c_4 D^2 n^r}
\,.
\]

Combining these, we find that except with probability $c_5 n! e^{-c_6 D^2
  n^r}$, we have simultaneously
\[  
\left| \frac{\Nind(H,\GG)}{I(n,k,r)} - 1 \right | < D
\qquad \text{ and } \qquad
\left| \frac{\Nord(H,\GG)}{I(n,k,r)} - \frac{1}{k!} \right | < D
\]
for all orderings $\ord_H$ and $\ord_V$. We now take $D := c_7
\sqrt{\frac{\log n}{n^{r-1}}}$ with $c_7$ chosen so that $c_5 n! e^{-c_6 D^2
  n^r} < 1$.
Then there is a $G$ satisfying the above bounds. The claim
then follows by the triangle inequality with $C:=2c_7$, since
\[
\left| \frac{\Nord(H,G)}{I(n,k,r)} - \frac{\Nord(H,G)}{\Nind(H,G)} \right |
= 
\frac{\Nord(H,G)}{\Nind(H,G)} \cdot \left| \frac{\Nind(H,G)}{I(n,k,r)} - 1
\right | < 1\cdot D\,.
\tag*{$\dashv$}
\]

This method of proof can be applied to many other classes
of structures, thereby establishing the uniqueness of consistent random
(admissible) orderings for these classes. These include: (i) the
Fra\"{i}ss\'e class of finite tournaments and (ii) the Fra\"{i}ss\'e class
of arbitrary $L$-structures for any finite language $L$ containing only
relation symbols of arity $\geq 2$. In both these cases, the uniform
ordering is the unique consistent random ordering. For another example,
consider a finite language $L$ containing at least one relation symbol of
arity $\geq 2$ and unary relation symbols $P_1,\dots , P_k$ and consider
the Fra\"{i}ss\'e class of structures for this language in which the
$P_1,\dots , P_k$ form a partition. In this case, the admissible orderings
for such a structure $\bfa$ will turn out to be those $<$ for which
$P_1^\bfa <  \dots  < P_k^\bfa$ and again the uniform ordering is the
unique consistent random admissible ordering.  (This also holds if the
language contains no relation symbols of arity $\geq 2$, but uniqueness is
straightforward in this case and does not need the methods of this section;
see \rref s.ot/.) Similarly, take as $\k$ the Fra\"{i}ss\'e class of
bipartite graphs with distinguished parts, i.e., structures of the form
$\langle V, E, A_1, A_2\rangle$, where $\langle V, E\rangle$ is a graph and
$A_1, A_2$ is a partition of $V$ such that there are no edges within each
$A_i$. Then for each such graph the relevant collection of admissible
orderings consists of all orderings $<$ such that $A_1 < A_2$ and the
uniform ordering is the unique consistent random admissible ordering.

\section{Bounded Degree Graphs (and Hypergraphs?)} \label{s.lower}

In this section, we construct non-uniform consistent random orderings on
graphs with bounded degrees, and estimate their total variation distance
from uniform. We believe we have a construction for hypergraphs as well,
but lack a proof.

\begin{thm}\label{T.approx_lower}
  There is a constant $C > 0$ with the following property.
  Let $n \ge 2$ be an integer.
  There is a consistent random ordering $(\mu_G)_G$ defined on all
  graphs of size at most $n$ such that 
  for every $k \in [3, n]$, there exists a
  graph $H$ of size $k$ with
  $d_\TV(\mu_H,\nu_H) \ge C/n$.
\end{thm}

In fact, we prove the following more general lower bound.

\procl t.hyper-lower
  There is a constant $C > 0$ with the following property.
  Let $D \ge 2$ be an integer.
  There exists a consistent random ordering $(\mu_G)_G$ defined on all
  graphs of degree at most $D$ such that for every $k \ge
  3$, there exists an 
  graph $H$ of degree at most $D$ and of size $k$
  such that
  $d_\TV(\mu_H,\nu_H) \ge C/D$.
\endprocl

\proof
Let $G=\langle V,E \rangle$ be a graph with maximal degree $D$. Let us make
every vertex $x\in V$ have degree exactly $D$ by adding $D-d_x$ additional
edges connecting $x$ to new, auxiliary vertices. Call the resulting graph
$\widehat{G}$.  Let $Z(e)$ be independent standard normal random variables for
the edges $e$ of $\widehat{G}$.  Define $Y(x) := \sum_{e \ni x} Z(e)$ for
vertices $x$ of $G$. (We do not bother defining $Y$ for vertices of
$\widehat{G} \setminus G$.)  Assign to the vertices of $G$ the order induced from
$Y \subset \bbR$.

If $H$ is an induced subgraph of $G$, then the inclusion of $H$ in $G$
can be extended to
a map from $\widehat{H}$ to $\widehat{G}$ that is 1-1 on edges (though some of
the vertices added to $H$ may be mapped to the same vertex of $G$.) Now the
IID Gaussians associated with the edges of $\widehat{G}$ can be pulled back to
$\widehat{H}$. This gives the ordering of $H$ as the restriction of the
ordering of $G$, thus showing that this ordering is consistent.

However, this ordering is not uniform. Given $k$, let $H$ be the graph on
vertices $\{1,\dots,k\}$ with only two edges, $e_1:=(1,2)$ and $e_2:=(2,3)$.
Let $A$ be the event that $Y(1) < Y(2) < Y(3)$.  Then $\nu_H(A) = 1/6$,
whereas $\mu_H(A) \geq 1/6 + 1/(6D) + o(1/D)$.
To see this latter fact, define $W_i := Y(i) - Z(e_1) - Z(e_2)$, so that
$A = \{W_1 < W_{2} < W_3\}$. Note that $W_1$, $W_2$, and $W_3$ are
independent normal random variables with variances $D$, $D-2$, and $D$,
respectively. Therefore, $\P(A) = \P[Z_1 < (1-2/D) Z_2 < Z_3]$ for 
independent standard normal random variables $Z_1$, $Z_2$, and $Z_3$.
Define $f(\epsilon) := \P[Z_1 < (1-\epsilon)Z_2 < Z_3]$. It suffices to
show that $f'(0) = 1/(2\pi\sqrt{3}) > 1/12$.

Write $\varphi(x)$ for the standard normal probability density.
Then for $a<b$,
\[
\frac{d}{d\epsilon} \P[a < (1-\epsilon)Z < b]\Big|_{\epsilon = 0}
= b\varphi(b) -
a\varphi(a),
\]
whereas for $a>b$ the derivative is trivially 0. Thus
\begin{align*}
\frac{d}{d\epsilon} \P[Z_1 < (1-\epsilon)Z_2 < Z_3]\Big|_{\epsilon = 0} 
&= \iint_{z_3>z_1} \big(z_3 \varphi(z_3) - z_1 \varphi(z_1)\big)
\varphi(z_1) \varphi(z_3) \,dz_1 \,dz_3 \\
&= \frac{1}{2\pi\sqrt{3}}\,.
\tag*{$\dashv$}
\end{align*}

An even easier construction holds for orderings of finite {\em connected}
bipartite graphs (such as trees) whose parts are not distinguished.
There we have the following:

\procl t.bipartite-non-uniform
There is a non-uniform consistent random ordering of the class of finite 
connected bipartite graphs.
\endprocl

\proof
We essentially saw this at the end of \rref s.quant/:
Given a finite connected bipartite graph,
let $A$ and $B$ be its two parts, named in random order. Order all of $A$
uniformly and all of $B$ uniformly independently, making all of $A$ less
than all of $B$.
It is easy to check that this is consistent. 
\Qed

\begin{question}
  Is there a non-uniform consistent random ordering on all finite bipartite
  graphs? Is there a non-uniform consistent random ordering on finite forests?
  Is there a non-uniform consistent random ordering on finite planar graphs?
\end{question}

\procl r.non-uniform-hyper
  We believe that the following extension of Theorem \ref{T.approx_lower} to
  uniform hypergraphs holds.  Let $n > r \ge 2$ be integers.
  Then there is a constant $C(r) > 0$ 
  and a consistent random ordering $(\mu_G)_G$ defined on all
  $r$-uniform hypergraphs of size at most $n$ such that 
  for every $k \in [r+1, n]$, there exists an 
  $r$-uniform hypergraph $H$ of size $k$ with
  $d_\TV(\mu_H,\nu_H) \ge C(r)/n^{r-1}$.
In fact, we believe the following more general lower bound.
Note that the \dfnterm{degree} of a vertex in a hypergraph is defined to be
the number of hyperedges that contain the vertex.
  Let $r, D \ge 2$ be integers.
  Then there is a constant $C(r) > 0$ with the following property:
  There exists a consistent random ordering $(\mu_G)_G$ defined on all
  $r$-uniform hypergraphs of degree at most $D$ such that for every $k \ge
  r+1$, there exists an 
  $r$-uniform hypergraph $H$ of degree at most $D$ and of size $k$
  such that
  $d_\TV(\mu_H,\nu_H) \ge C(r)/D$.
\endprocl

It appears via simulations that the following modification of the proof of
Theorem \ref{T.approx_lower} should work.
The consistency condition (ii) means that for $r \ge 3$, assigning a single
Gaussian variable to each hyperedge and summing the variables of the edges
containing a vertex does {\it not} yield a consistent ordering. The following
construction {\it is} consistent and appears to give the claimed lower bound.
Let $Z := \Seq{Z_1, Z_2, \ldots, Z_r}$ be a collection of $(r-1)$-independent
exchangeable standard normal random variables that are not independent.
(Recall that $Z$ is \dfnterm{exchangeable} means that the law of $Z$ is
invariant under permutations of its coordinates.) For
example, let $\Phi$ be the standard normal cumulative distribution
function. Let $T_1, \ldots, T_r$
be IID uniform $[0, 1]$ random variables conditioned to sum to 0 mod 1. Now
define $Z_i := \Phi^{-1}(T_i)$.
Let $G$ be an $r$-uniform hypergraph with maximal degree $D$.  Let $Z(e) =
\Seq{Z(e)_x \st x \in e}$ be
independent copies of $Z$ for the hyperedges $e$ of $G$. Also,
for each vertex $x$ of $G$ of degree $d_x < D$, let $Z(x, 1), \ldots,
Z(x, D-d_x)$ be additional independent standard normal random variables.
Define $Y(x) := \sum_{e \ni x} Z(e)_x + \sum_{i=1}^{D - d_x} Z(x, i)$ for
vertices $x$ of $G$.  Note that the process $Y$ is $\Aut(G)$-invariant.
Assign the vertices of $G$ the order induced from $Y
\subset \bbR$.  This is consistent and appears not to be uniform.

In fact, given $k$, let $H$ be the hypergraph consisting of the vertices
$\{x_1, x_2, \ldots, x_k\}$ and two hyperedges, $e_1 := \{x_1, x_2, \ldots,
x_r\}$ and $e_2 := \{x_2, x_3, \ldots, x_{r+1}\}$.
Let $A$ be the event that $Y(x_1) < Y(x_2) < \cdots < Y(x_r) < Y(x_{r+1})$.
Then $\nu_H(A) = 1/(r+1)!$, whereas it seems that
$\mu_H(A) \ge 1/(r+1)! + C(r)/D + o(1/D)$.

\bsection{Dense Hypergraphs of Large Girth}{s.G0}

Here we prove a version of a lemma of \rref b.NR/ that will be very
useful to us in analyzing more complicated structures in the next two
sections. 
There are various kinds of paths one can define in a hypergraph.
We use the following.
A \dfnterm{path}
in a hypergraph is an alternating sequence $\Seq{x_1, e_1, x_2, e_2,
\ldots, x_L, e_L, x_{L+1}}$ of vertices $x_i$ and
hyperedges $e_i$ such that $x_i \ne x_{i+1}$,
$x_i, x_{i+1} \in e_i$, and $e_i \ne e_{i+1}$ for
all $i \in [1, L]$. Such a path is said to \dfnterm{join} $x_1$ to $x_{L+1}$,
to have \dfnterm{length} $L$, and to be a
\dfnterm{cycle} if $L \ge 2$ and $x_1 = x_{L+1}$.
The \dfnterm{girth} of a hypergraph is the minimal length
of a cycle that it contains.
A hypergraph is \dfnterm{connected}
if every pair of distinct vertices is joined by
some path.

\procl l.G_0
Let $r \ge 2$ and $g \ge 3$ be integers. There is a constant $C = C(r, g)$
so that for all $n \ge r$, there
exists an $r$-uniform hypergraph on $n$ vertices and at least
$C n^{(g-1)/(g-2)}$
hyperedges that has girth at least $g$.
\endprocl

\proof
The (standard) method is to take a random hypergraph, and remove all edges
that are in short cycles. Let $c_i$ denote constants that depend on $r$
and $g$, but not on $n$.
Let $p := a/n^{r-(g-1)/(g-2)}$ for a small constant $a<1$ to be chosen later.
Let $\GG$ be the random $r$-uniform hypergraph on $n$ vertices such that each
hyperedge belongs to $\GG$ independently with probability $p$.
Thus, the expected number of hyperedges in $\GG$ is ${n \choose r} p$,
which is at least $c_1 a n^{(g-1)/(g-2)}$.
Let $2 \le j < g$.
The union of the hyperedges of any minimal cycle of length $j$ contains
at most $r j - j$ vertices.
The number of cycles of length $j$ whose union is a given set of
size $i$ is at most ${i \choose r}^j$, and each such cycle has probability
$p^j$ to belong to $\GG$.
Also, the number of hyperedges that belong to some minimal cycle of length
$j$ is at most $j$ times the number of such minimal cycles.
Thus, the
expected number of hyperedges that belong to some minimal cycle of length
$j$ is at most $c_2 \sum_{i=r}^{rj-j} {n \choose i} p^j \le
c_3 n^{(r-1)j} p^j = c_3 a^j n^{j/(g-2)}$.
Hence the 
expected number of hyperedges that belong to some cycle of length
less than $g$ is at most $\sum_{j=2}^{g-1} c_3 a^j n^{j/(g-2)}
\le c_4 a^2 n^{(g-1)/(g-2)}$.
Now for $a$ sufficiently small, $C := c_1 a - c_4 a^2 > 0$. 
That is, for $a$ sufficiently small, the expected number of
hyperedges in $\GG$ that do not belong to any cycle of length less than
$g$ is more than $C n^{(g-1)/(g-2)}$.
Therefore, there is some hypergraph with more than $C n^{(g-1)/(g-2)}$
hyperedges that do not belong to any cycle of length less than $g$.
Take such a hypergraph and remove all hyperedges in cycles of length less than
$g$.
\Qed

We remark that the hypergraph may be constructed to be connected at the
price of allowing the number of vertices to be in the interval $[n,
n+r-2]$.
To see this, if the result above
is disconnected and has at least $r$ connected components, then we may add
a hyperedge to reduce the number of components without creating any new
cycles.
If, on the other hand, the number $p$ of connected components is between 2
and $r-1$, then we may add a hyperedge containing $r-p$ new vertices to
make it connected without creating any new cycles.

\bsection{Forbidden Subgraphs}{s.forb}

Given the edge set $E$ of a graph $K$, identify subsets $A \subseteq E$
with their indicator functions $\II A \in (\bbZ_2)^E$, so that $A_1 \symdiff
A_2$ is identified with $\II {A_1} + \II {A_2}$.
We say that a simple cycle $C$ is \dfnterm{generated} by simple
cycles $C_1, \ldots, C_j$ if $C$ is the sum (in the previous sense) of the
$C_i$ ($1 \le i \le j$), where we regard a simple cycle as its set of
edges.
Given an integer $g \ge 3$, say that a graph $K$ is \dfnterm{$g$-small} if
$K$ is connected, has no cutpoints, and all simple cycles in $K$ are
generated by simple cycles in $K$ of length $< g$. 
For example, if $K$ is connected, has no cutpoints, and has size $< g$,
then $K$ is $g$-small.
For another example, note that the usual Cayley graph of $\bbZ^2$, i.e., the
infinite square lattice graph, is 5-small.

Given a class $\forb$ of graphs,
write $Forb(\forb)$ for the class of finite graphs that have no induced
subgraph in $\forb$.
Note that if $\forb$ is finite and consists of connected finite graphs
without cutpoints, then $\forb$ contains only $g$-small graphs for some
fixed $g$.
Also, if $\forb$ is hereditary and each graph $K$ in $\forb$ has the
property that all simple cycles in $K$ are
generated by simple cycles in $K$ of length $< g$, then $Forb(\forb)
= Forb(\forb')$ for some class $\forb'$ that
contains only $g$-small graphs. Indeed, we may let $\forb'$ be the class of
graphs in $\forb$ that are connected and have no cutpoints.

\procl t.forb-uniform
  Let $g \ge 3$ be an integer and $\forb$ be a collection of $g$-small
  graphs.
  The uniform ordering is the unique consistent random ordering on the class
  $Forb(\forb)$.
\endprocl

The quantitative version of this theorem follows.
In it, we speak of a restricted class of (induced) embeddings of a graph $H$ in a
graph $G$.
We use the superscript {res} to denote the restriction in counting
embeddings and in counting ordered embeddings.
The restriction depends on both $H$ and $G$ and can be arbitrary, but it
does not depend on orderings of $H$ and $G$. We denote by $\Nindr(H, G)$
the number of restricted embeddings, while for any fixed orderings $<_H$, $<_G$ of
$H$, $G$, resp., we denote by $\Nordr(H, G)$ the number of restricted
embeddings that preserve $<_H$, $<_G$.
The proof that \rref t.forb/ implies \rref t.forb-uniform/ is the same as
that of \rref l.TVbd/. 

\procl t.forb
Let $k, g \ge 3$ be integers and $\forb$ be a collection of $g$-small
graphs.
There exists a constant $C(k, g)$ with the following property.
For every graph $H \in Forb(\forb)$ on $k$ vertices and $n \ge k$, 
there exists a graph $G \in Forb(\forb)$
on $n$ vertices such that there is a restricted class of embeddings for
which
$\Nindr(H, G) > 0$
and for every pair of orderings
$\ord_G$ of $G$ and $\ord_H$ of $H$, 
$$
\left|{\Nordr(H, G) \over \Nindr(H, G)} - {1 \over k!}\right|
<
\sqrt{\log n \over C(k, g) n^{1/(g-2)}}
\,.
$$
\endprocl

\proof
Let $L$ be the number of automorphisms of $H$.
Let $C(k, g)$ be the constant of \rref l.G_0/.
Write $m := \ceil{C(k, g) n^{(g-1)/(g-2)}}$.
Let $G_0$ be a $k$-uniform hypergraph on $n$ vertices $V$ having $m$
hyperedges $e_1, \ldots, e_m$ and girth at least $g$.
Such a $G_0$ exists by \rref l.G_0/ when $n \ge k$.
Since the girth of $G_0$ is larger than 2, no two hyperedges share more
than one vertex of $G_0$.
Let $\GG$ be the random graph obtained from $G_0$ as follows.
For $1 \le i \le m$, let $\HH_i$ be the random graph isomorphic to $H$ with
vertex set $e_i$ induced by a uniform random bijection between
$e_i$ and $\verts(H)$.
Choose $\HH_i$ independently.
The union of all $\HH_i$ is $\GG$.

We claim that $\GG \in Forb(\forb)$.
Suppose not and that $\GG$ has an induced subgraph $K \in \forb$.
Since $H \in Forb(\forb)$, it follows that $K$ is not contained entirely
within any $\HH_i$.
Since $K$ is connected and has no cutpoints, it also follows that there is
a simple cycle $C$ of $K$ that has length at least $g$ and that there
is a cycle $C' = \Seq{x_1, e_{r_1}, x_2, e_{r_2}, \ldots, x_L, e_{r_L},
x_{L+1}}$ in $G_0$ and paths $P_i \subseteq \HH_{r_i}$ joining $x_i$
to $x_{i+1}$ so that $C$ is the concatenation of the paths $P_1, \ldots,
P_L$.
Choose such a pair of cycles $C$, $C'$ with $L$ a minimum. Then all $x_i$
are distinct.
Now $C$ is generated by simple cycles in $K$ of length $< g$, each of
which, being simple, lies within some $\HH_i$.
However, the sum of cycles is an even graph, i.e., all its vertices have
even degree, whereas the intersection of $C$ with each $\HH_i$ is a union of
paths with distinct endpoints and thus is not an even graph.
Since this intersection must be generated by the cycles that lie within
$\HH_i$, we obtain a contradiction, which establishes our claim.

We shall restrict to the embeddings of $H$ in $\GG$ that embed $H$ in some
hyperedge of $G_0$. Thus, $\Nindr(H, \GG) = L m$.

Fix $\ord_V$ and $\ord_H$. Let $G$ be a possible value of $\GG$.
Let $ f(G) := \Nordr(H, G)/(L m)$.
Now $\Nordr(H, \GG)$ has a binomial distribution with parameters $(m, L/k!)$.
Thus, Chernoff's inequality yields
$$
\P\big[|f(\GG) - 1/k!| \ge D\big]
\le
2 \exp\{-2L^2 D^2 m\}
\,.
$$
Choose 
$$
D := \sqrt{\frac{n \log n}{m}} \le \sqrt{\frac{\log n}{C(k, g) n^{1/(g-2)}}}
\,.
$$
This gives
$$
\P\big[|f(\GG) - 1/k!| \ge D\big]
\le
2 \exp\{-2 L^2 n \log n\}
\,.
$$
Since this holds for every $\ord_V$ and every $\ord_H$
and the number of pairs of
orderings of $\verts(G_0)$ and $\verts(H)$ is $n! k!$, we obtain the result.
\Qed

A similar proof clearly works for the Fra\"{i}ss\'e class of $\a$-free hypergraphs of a given type and also for the Fra\"{i}ss\'e class of abstract simplicial complexes of bounded dimension. So we have:

\procl t.free
The uniform ordering is the unique consistent random ordering on the class
of finite $\a$-free hypergraphs of a given type and for the class of finite
simplicial complexes of bounded dimension.
\endprocl

\procl r.methods The method used in this section to prove uniqueness of
consistent random orderings works also for the classes considered in \rref
s.quant/.
\endprocl

\bsection{Metric Spaces}{s.metric}

Since there are many interesting classes of metric spaces, they provide a
fertile ground for investigation. First we prove that when (essentially) no
restriction is placed on the class of metric spaces, the only consistent random
ordering is the uniform one. In the next section,
we show that some particular classes have
other consistent orderings.

\procl t.metric-uniform
Let $S$ be an additive subsemigroup of\/ $\bbR^+$.
The uniform ordering is the unique consistent random ordering on the class
of metric spaces with non-zero distances in $S$.
\endprocl

A quantitative version follows.

\procl t.metric
Let $k, \alpha \ge 3$ be integers.
There exists a constant $C'(k, \alpha)$ with the following property.
For every metric space $\langle X, d\rangle$ on $k$ vertices satisfying $d(x_1, x_2) \le
\alpha d(x_3, x_4)$ for all $x_i \in X$ with $x_3 \not= x_4$, and for every $n \ge k$,
there exists a metric space $Y$
on $n$ vertices with (non-zero) distances in the additive semigroup generated by the (non-zero) distances
in $X$ and such that 
there is a restricted class of embeddings for which
$\Nindr(X, Y) > 0$
and for every pair of orderings
$\ord_X$ of $X$ and $\ord_Y$ of $Y$, 
$$
\left|{\Nordr(X, Y) \over \Nindr(X, Y)} - {1 \over k!}\right|
\le
\sqrt{\log n \over C'(k, \alpha) n^{1/{(\alpha-1)}}}
\,.
$$
\endprocl

\proof
Let $L$ be the number of isometries of $X$, so $L \le k!$.

Let $C'(k, \alpha) := C(k, \alpha+1)$, where $C(r, g)$ is
the constant of \rref l.G_0/.
Write $m := \ceil{C'(k, \alpha) n^{\alpha/(\alpha-1)}}$.
Let $G_0$ be a $k$-uniform hypergraph on $n$ vertices $V$ having $m$
hyperedges $e_1, \ldots, e_m$ and girth at least $\alpha+1$.
Such a $G_0$ exists by \rref l.G_0/ when $n \ge k$.
Note in particular that no two hyperedges share more than one vertex of
$G_0$.
Let $\YY$ be the random metric space on $V$ obtained from $G_0$ as follows.
For $1 \le i \le m$, let $(\XX_i, d_i)$ be the random metric space isometric to
$X$ on $e_i$ induced by a uniform random bijection between
$e_i$ and $X$.
Extend the resulting metric $d$ to all pairs
$z, w \in \verts$ by using the induced shortest-path metric for those
pairs that are joined by a path of points whose consecutive distances have
already been defined.
Let $\beta$ be the maximum distance thereby obtained and
define $d(z, w)$ to be $\beta$ if there is no
path of $G_0$ that joins $z$ and $w$.
This defines $\YY$.
Note that the restriction of $d$ to each $\XX_i$ agrees with $d_i$ since if
$P$ is a shortest path between two points $x, y \in \XX_i$, then either
$P$ is contained within $\XX_i$, in which case the result follows from the
triangle inequality for $d_i$, or $P$ has length at least $\alpha$, in
which case the length of $P$ is at least $\alpha d_i(x, y)/\alpha = d_i(x,
y)$ by the hypothesized inequality involving $\alpha$.

We restrict to the embeddings of 
$X$ in $\YY$ that embed $X$ in some hyperedge. Thus,
$\Nindr(X, \YY) = L m$.

Fix $\ord_Y$ and $\ord_X$. Let $Y$ be a possible value of $\YY$.
Let $f(Y) := \Nordr(X, Y)/(L m)$. 
Now $\Nordr(X, \YY)$ has a binomial distribution with parameters $(m, L/k!)$.
Thus, Chernoff's inequality yields
$$
\P\big[|f(\YY) - 1/k!| \ge D\big]
\le
2 \exp\{-2L^2 D^2 m\}
\,.
$$
Choose $D := \sqrt{n \log n/m}$.
This gives
$$
\P\big[|f(\YY) - 1/k!| \ge D\big]
\le
2 \exp\{-2 L^2 n \log n\}
\,.
$$
Since this holds for every $\ord_Y = \ord_{\verts}$ and every $\ord_X$
and the number of pairs of
orderings of $\verts$ and $\ord_X$ is $n! k!$, we obtain the result.
\Qed


\bsection{Euclidean and Other Metric Spaces}{s.eucl}

For certain classes of metric spaces, there is a non-uniform consistent
random ordering.
We begin by describing an idea of Leonard Schulman (personal communication,
2010) for randomly ordering finite subsets $X$ of Euclidean space $\bbR^n$.
Project $\bbR^n$ orthogonally onto a uniformly
random oriented line through the origin. Order the points of $X$
corresponding to the order of their projections
on the line. Write $\mu_X$ for the law
of this random order.
By considering translations, rotations, and reflections separately,
it is not hard to see that if $\phi : \bbR^n \to \bbR^n$ is an
isometry, then $\mu_{\phi[X]} = \phi_* (\mu_X)$; furthermore,
$\mu_Y$ is induced by
restriction of $\mu_X$ when $Y \subseteq X$. In that sense,
$X \mapsto \mu_X$ is consistent for the class of finite subsets of $\bbR^n$.
In addition, it is consistent in the following sense:
if $X \subset \bbR^m \subset \bbR^n$ for some $m < n$, then the probability
measure on orderings of $X$ induced by random projections of $\bbR^m$ is
the same as the one induced from $\bbR^n$. 
Finally, it is consistent in that whenever $\phi : X \to Y$ is an isometry
of finite subsets of $\bbR^n$ with their induced metrics, then $\mu_Y =
\phi_*(\mu_X)$ since $\phi$ extends to an isometry of all of $\bbR^n$.
(That is, Euclidean space is ultrahomogeneous.)

To extend this idea, call a metric space \dfnterm{Euclidean} if it is
isometric to a subset of some Euclidean space.

\procl t.Eucl-consis
Fix an injection $f : [0, \infty) \to [0, \infty)$ with $f(0) = 0$.
There is a non-uniform consistent random ordering on the class of
those finite metric spaces $\langle X, d\rangle$ for which
$\Seq{X, f \circ d}$ is Euclidean.
\endprocl

\proof
Given an isometric embedding $\phi : \Seq{X, f \circ d} \to \bbR^{|X|}$,
define the ordering $\overline\mu_X := \phi^* \mu_{\phi[X]}$ pulled back
from the ordering on the image of $X$ defined above.
By ultrahomogeneity, $\overline\mu_X$ is independent of choice of $\phi$.
However, $\overline\mu_X$
is not uniform whenever not all non-0 distances in $X$ are the same. 
\Qed

Examples include metric spaces of \dfnterm{negative type}, which can be
defined as those metric spaces $\Seq{X, d}$ such that $\Seq{X, f \circ d}$
is Euclidean for $f(s) := \sqrt s$. In fact, it then turns out that one may
also take $f(s) := s^\alpha$ for any $\alpha\in (0, 1/2]$; see \cite{Sch}.
Examples of metric spaces of negative type include ultrametric
spaces, spheres, hyperbolic spaces, and all $L^p$ spaces for $1 \le p \le
2$; see \cite[Theorem 3.6]{Mec} and the references there.

It should be pointed out that the Ramsey properties of classes of Euclidean metric spaces are far from being understood and it is conceivable that the above have some relevance in this context.

One can also establish \rref t.hyper-lower/,
i.e., that there is a non-uniform consistent random ordering on graphs of bounded
degree, by this random projection method. 
In order to choose distances on the vertices of a graph $G = \langle V, E\rangle$ with
maximum degree at most $D$ that give a metric
of negative type, fix $a < 1$ such that $a/(1-a) = D$ and for $x, y \in
V$, define
\[
d(x, y) := \begin{cases}
            0 &\mbox{if } x = y,\\
            1 &\mbox{if } \{x, y\} \in E,\\
            a &\mbox{otherwise}.\\
            \end{cases}
\]
We use the fact \cite{Sch} that a finite metric space has negative type iff
its matrix of distances is \dfnterm{conditionally negative semidefinite
(CNSD)}. 
Here, a matrix $M$ is called CNSD if $(M v, v) \le 0$ for all vectors $v$
that are orthogonal to the constant vectors.
Let $A$ be the adjacency matrix and $J$ be the all-ones matrix indexed by the
vertices. Then the distance matrix is $A + a (J-A-I) = (1-a)A - a
I + a J$. When this matrix acts on a vector orthogonal to the
constants, $J$ sends the vector to 0, so the distance matrix is CNSD iff
$(1-a)A - a I$ is CNSD iff $-(D I - A)$ is CNSD. But in
fact, this matrix is \dfnterm{negative semidefinite (NSD)} as seen, e.g., by
comparison to the graph Laplacian matrix.
Here, a matrix $M$ is called NSD if $(M v, v) \le 0$ for all vectors $v$.

It seems that using random distances, one can also use this method for
uniform hypergraphs, but this appears even harder than the method used in
the discussion of \rref r.non-uniform-hyper/.

If one wishes, one can define the random ordering on finite subsets of
Euclidean spaces in all dimensions at once by considering finite
subsets $X$ of Hilbert space $\Hilb := \ell^2(\bbN)$ instead. For that case,
let $\Seq{v_n \st n \ge 0}$ be any orthonormal basis of $\Hilb$ and let
$\Seq{Z_n \st n \ge 0}$ be independent standard normal random variables.
Order $X$ by the order on the real numbers $\sum_n Z_n \langle x, v_n
\rangle$ for $x \in X$; this sum converges a.s.\ by Kolmogorov's
Three-Series Theorem because $\sum_n |\langle
x, v_n \rangle|^2 < \infty$.
The sum has a normal distribution with variance $\|x\|^2$.
The spherical symmetry of the standard multivariate normal distribution in
Euclidean spaces shows that this random ordering does not depend on the
choice of basis and agrees with the random ordering $\mu_X$ previously defined.

This is related to the L\'evy-Ciesielski construction of Brownian motion:
First, note that $\bbR^+$ has negative type, as we can see by
embedding $\bbR^+$ into $L^2(\bbR^+)$ via 
mapping $x \in \bbR^+$ to the function
$s \mapsto \bfone_{[0, x]}(s)$. Second, identify $L^2(\bbR^+)$ with $\Hilb$ by
taking as orthonormal basis the Haar basis. Then it is not hard to see that
$\sum_n Z_n \langle \bfone_{[0, x]}, v_n \rangle$ a.s.\ converges uniformly
for $x$ belonging to any compact set. Since it is a mean-0 Gaussian process
with covariance $(x, y) \mapsto \langle \bfone_{[0, x]}, \bfone_{[0, y]}
\rangle = \min\{x, y\}$, it is standard Brownian motion.
This is precisely the L\'evy-Ciesielski construction.

This concludes the discussion of consistent random orderings on specific classes of finite structures. The rest of the paper will connect this with the unique ergodicity phenomena.

\section{Topological Dynamics and Unique Ergodicity}

We prove here some facts concerning unique ergodicity for amenable groups that will be used in subsequent sections.
First we note the following general fact.

\procl p.Char
Let\/ $\gp$ be a topological group. Then $\gp$ is amenable iff the universal minimal flow of $\gp$ admits an invariant measure. Moreover, $\gp$ is uniquely ergodic iff the universal minimal flow of $\gp$ is uniquely ergodic.
\endprocl

\proof The first statement is obvious, since every flow contains a minimal
subflow. For the second, it is enough to show that if $X$ is a uniquely
ergodic $\gp$-flow, $Y$ is a $\gp$-flow and $\pi : X \rightarrow Y$ is a surjective, continuous $\gp$-map, then $Y$ is uniquely ergodic. For that it is enough again to show that if $\nu$ is a $\gp$-invariant measure on $Y$, then there is a $\gp$-invariant measure $\hat{\nu}$ such that $\pi_* \hat{\nu} = \nu$.

First note that there is a measure $\mu_0$ on $X$ such that $\pi_*\mu_0 =
\nu.$ Indeed the set 
$$
\{\pi_*\mu \st \mu \text{ is a measure on } X \}
$$
is a compact, convex set, containing all Dirac measures, thus, by
Hahn-Banach, it contains all measures on $Y$. It follows that the set of
all measures $\mu$ on $X$ with $\pi_* \mu = \nu$ is a non-empty, compact,
convex set of measures on which $\gp$ acts continuously by affine
transformations (the action given as usual by: $\int f d (g\cdot \mu)  =
\int (g^{-1}\cdot f) \,d\mu$, for $g\in \gp$, $f\in C(X)$, where 
$g\cdot f (x) = f(g^{-1}\cdot x)$). Since $\gp$ is amenable, this action
has a fixed point $\hat{\nu}$ (see \cite[G.1.7]{BHV}) and thus $\hat{\nu}$ is as required.
\Qed

\medskip
Next we provide a characterization of unique ergodicity in the case of compactly approximable groups.

Let $\gp$ be a compactly approximable topological group and let $\Delta_0\subseteq
\Delta_1\subseteq\cdots$ be a sequence of compact subgroups with $\overline
{\bigcup_n\Delta_n}=\gp$.  Let $\mu_n$ be the Haar measure of $\Delta_n$.  Let
$X$ be a $\gp$-flow.  For $f\in C(X), n\in\bbN$, define
the \dfnterm{averaging operator} $A_n$ on $C(X)$ by
\[
A_n(f)(x):=\int_{\Delta_n}f(g^{-1}\cdot x)\,d\mu_n(g).
\]
The following is an analog of 4.9 in Glasner \cite{Gl}.

\procl t.aver
In the preceding notation, the following are equivalent:
\begin{itemize}
\item[(i)] $\forall f\in C(X)\ \exists f^*\in\bbC$ ($A_n(f)$ converges
uniformly to $f^*$),
\item[(ii)] $\forall f\in C(X)\ \exists f^*\in \bbC$ ($A_n(f)$ converges
pointwise to $f^*$),
\item[(iii)] There is a $\gp$-invariant measure $\mu$ on $X$ such that
\[
\forall f\in C(X)\ (A_n(f)\text{ converges pointwise to }\int f\,d\mu )
\]
\item[(iv)] The $\gp$-flow $X$ is uniquely ergodic.
\end{itemize}
\endprocl

\proof (i) $\Rightarrow$ (ii) is obvious.

\smallskip
(ii) $\Rightarrow$ (iii):  Put $A(f):=\lim_{n\to\infty}A_n(f)(x)\in\bbC$
(for any $x\in X$).  Then $A$ is a positive linear functional on
$C(X)$ with $A(1)=1$, so by the Riesz Representation Theorem
there is a measure $\mu$ on $X$ with
\[
A(f)=\int f\,d\mu .
\]
We shall show that $\mu$ is $\gp$-invariant.  It is of course enough to
show that it is ($\bigcup_n\Delta_n$)-invariant or equivalently that $A$ is
($\bigcup_n\Delta_n$)-invariant, where $\gp$ acts on $C(X)$ by $g\cdot
f(x):= f(g^{-1}\cdot x)$.  Fix $h\in\bigcup_n\Delta_n$.  Then
\begin{align*}
A(h\cdot f)(x) &=\lim_{n\to\infty}\int_{\Delta_n}(h\cdot f)(g^{-1}\cdot x)
\,d\mu_n(g)\\
&=\lim_{n\to\infty}\int_{\Delta_n}f(h^{-1}g^{-1}\cdot x)\,d\mu_n(g)\\
&=\lim_{n\to\infty}\int_{\Delta_n}f\big((gh)^{-1}\cdot x\big)\,d\mu_n(g)\\
&=\lim_{n\to\infty}\int_{\Delta_n}f(g^{-1}\cdot x)\,d\mu_n(g)\\
&=A(f)(x)
\end{align*}
by the invariance of Haar measure.

\smallskip
(iii) $\Rightarrow$ (iv): Let $\mu$ be a $\gp$-invariant measure with 
$A_n(f)(x)\to\int f\,d\mu =f^*$ pointwise.  Let $\nu$ be any $\gp$-invariant
measure.  We shall show that $\mu =\nu$.  By Lebesgue Dominated
Convergence, we have 
\[
\int A_n(f)(x)d\nu (x)\to\int f^*d\nu =f^*=\int f\,d\mu .
\]
But also by Fubini and the $\gp$-invariance of $\nu$,
\begin{align*}
\int A_n(f)(x)d\nu (x) &=\int \left (\int_{\Delta_n}f(g^{-1}\cdot x)\,d\mu_n
(g)\right )d\nu (x)\\
&=\int_{\Delta_n}\left (\int f(g^{-1}\cdot x)d\nu (x)\right )\,d\mu_n(g)\\
&=\int_{\Delta_n}\left (\int fd\nu \right )\,d\mu_n(g)=\int f\,d\nu ,
\end{align*}
i.e., $\int f\,d\mu =\int f\,d\nu$, $\forall f\in C(X)$, so $\mu
=\nu$.

\smallskip
(iv) $\Rightarrow$ (i):  Let $\mu$ be the unique $\gp$-invariant
measure. If (i) fails, there is $f\in C(X)$ such that $A_n(f)$ does
not converge uniformly to $\int f\,d\mu$.  Fix then $\epsilon >0$ such
that for infinitely many $n$, there is $x_n\in X$ with $|A_n(f)
(x_n)-\int f\,d\mu |\geq\epsilon$.  For such $n$, let $\rho_n$ be the 
measure on $X$ defined by
\[
\int h\,d\rho_n:=\int_{\Delta_n}h(g^{-1}\cdot x_n)\,d\mu_n(g)
\]
for $h\in C(X)$.  Thus
\[
\int f\,d\rho_n=A_n(f)(x_n),
\]
so
\[
\left |\int f\,d\rho_n-\int f\,d\mu\right |\geq\epsilon .
\]
By compactness there is a measure $\rho_\infty$
and a subsequence $(\rho_{n_i})$ converging to $\rho_\infty$ in the 
weak$^*$-topology of measures.  Thus
\[
\left |\int f\,d\rho_\infty -\int f\,d\mu\right |\geq\epsilon ,
\]
so $\rho_\infty\neq\mu$.  We shall show that $\rho_\infty$ is 
$\gp$-invariant, which is a contradiction.  Fix $g_0\in\bigcup_n\Delta_n$.
Then for any $h\in C(X)$,
\begin{align*}
\int (g_0\cdot h)d\rho_\infty &=\lim_{i\to\infty}\int (g_0\cdot h)
d\rho_{n_i}\\
&=\lim_{i\to\infty}\left [\int_{\Delta_{n_i}}(g_0\cdot h)(g^{-1}\cdot
x_{n_i})\,d\mu_{n_i}(g)\right ]\\
&=\lim_{i\to\infty}\left [\int_{\Delta_{n_i}}h\big((gg_0)^{-1}\cdot
x_{n_i}\big)
\,d\mu_{n_i} (g)\right ]\\
&=\lim_{i\to\infty}\int_{\Delta_{n_i}}h(g^{-1}\cdot x_{n_i})\,d\mu_{n_i}(g)\\
&=\lim_{i\to\infty}\int h\,d\rho_{n_i}
\\
&=\int h\,d\rho_\infty .
\tag*{$\dashv$}
\end{align*}

\bsection{Universal Minimal Flow of Automorphism Groups and Unique Ergodicity}{s.S2}

Throughout the rest of this paper we shall work in the following
context developed in \cite{KPT} (although our notation will be slightly
different).

Consider a Fra\"{i}ss\'e class $\k$ in a language $L$.  Let $L^*:=L
\cup\{<\}$ be the language obtained by adding a binary relation symbol 
$<$ to $L$.  A structure $\bfa^*$ for $L^*$ has the form $\bfa^*=
\langle\bfa ,<_\bfa\rangle$, where $\bfa$ is a structure for $L$ and $<_\bfa$
is a binary relation on $A$ (= the universe of $\bfa$).  We often write
more simply $\Seq{\bfa, <}$ for $\Seq{\bfa, <_\bfa}$. A class
$\k^*$ of finite structures in $L^*$ is called an \dfnterm{order class}
if ($\langle\bfa ,<\rangle\in\k^*\Rightarrow\ <$ is a linear ordering
on $A$).  For such $\bfa^*=\langle\bfa ,<\rangle$, let $\bfa^*\restrict L:=
\bfa$.

We say that an order class $\k^*$ on $L^*$, closed under isomorphism, is an
\dfnterm{order expansion} of $\k$ if $\k =\k^*\restrict
L:=\{\bfa^*\restrict L\st\bfa^*\in\k^*\}$.  In this case, if $\bfa\in\k$
and $\bfa^*:=\langle\bfa ,<\rangle\in\k^*$, we say that $<$ is a
$\k^*$-\dfnterm{admissible ordering} of $\bfa$.  The order expansion $\k^*$
of $\k$ is \dfnterm{reasonable} if for every $\bfa ,\bfb\in\k$ with
$\bfa\subseteq\bfb$ and any $\k^*$-admissible ordering $<_\bfa$ on $\bfa$,
there is a $\k^*$-admissible ordering $<_\bfb$ on $\bfb$ such that
$<_\bfa\subseteq <_\bfb$.

If $\k$ is a Fra\"{i}ss\'e class with $\bfk ={\rm Flim}(\k )$ and $\k^*$
is a reasonable order expansion of $\k$, we denote by $X_{\k^*}$ the 
space of linear orderings $<$ on $K$ such
that for any finite substructure $\bfa$ of $\bfk$, $<\restrict A$ is 
$\k^*$-admissible on $\bfa$.  We call these the $\k^*$-\dfnterm{admissible
orderings} on $\bfk$.  They form a compact, non-empty subspace of
$2^{K^2}$, with the product topology, on which the group $\gp:={\rm
Aut}(\bfk )$ acts continuously. Thus $X_{\k^*}$ is a $\gp$-flow.

If $\k^*$ is an order expansion of $\k$, we say that $\k^*$ satisfies
the \dfnterm{ordering property (OP)} if for every $\bfa\in\k$, there is
$\bfb\in\k$ such that for every pair of $\k^*$-admissible orderings $<_\bfa$ on
$\bfa$ and $<_\bfb$ on $\bfb$, $\langle\bfa ,<_\bfa\rangle$ can be embedded
in $\langle\bfb ,<_\bfb\rangle$.

We also say that a class $\k$ of finite structures in a given language
$L$ has the \dfnterm{Ramsey Property (RP)} if for any $\bfa ,\bfb\in\k$
with $\bfa\leq\bfb$ (i.e., $\bfa$ can be embedded in $\bfb$), there
is $\bfc\in \k$ with $\bfb\leq\bfc$ such that for any coloring $c\colon
\begin{pmatrix}\bfc \\ \bfa \end{pmatrix}\to \{1,\dots ,k\}$, there is
$\bfb'\in\begin{pmatrix}\bfc \\ \bfb\end{pmatrix}$ such that 
$c\restrict \begin{pmatrix}\bfb' \\ \bfa\end{pmatrix}$ is constant.  Here
for $\bfd \leq\bfe$, we let $\begin{pmatrix}\bfe \\
\bfd\end{pmatrix}:=$ the set of all substructures of $\bfe$
isomorphic to $\bfd$.

The following is shown in \cite[7.4, 10.8]{KPT}.
\procl t.KPT \cite{KPT}
Let $\k$ be a Fra\"{i}ss\'e class in $L$ and $\k^*$ a 
reasonable order expansion of $\k$ in $L^*$ that is also a 
Fra\"{i}ss\'e class.  Let $\bfk :={\rm Flim}(\k )$ and $\gp:={\rm
Aut}(\bfk )$. Then the following are equivalent:
\begin{enumerate}
\item[(i)] $X_{\k^*}$ is a minimal $\gp$-flow;
\item[(ii)] $\k^*$ satisfies the {\rm OP}.
\end{enumerate}
Moreover, the following are equivalent:
\begin{enumerate}
\item[(a)] $X_{\k^*}$ is the universal minimal flow of\/ $\gp$;
\item[(b)] $\k^*$ satisfies the {\rm OP} and the {\rm RP}.
\end{enumerate}
\endprocl

We call $(\k ,\k^*)$ an \dfnterm{excellent pair} if $\k$ is a
Fra\"{i}ss\'e class in $L$, $\k^*$ is a  Fra\"{i}ss\'e, reasonable order
expansion of $\k$ in $L^*$, and $\k^*$ satisfies the {\rm OP} and the {\rm
RP}. It was shown in \cite[9.2]{KPT} that  if $\k$ admits an excellent
pair $(\k ,\k^*)$, then $\k^*$ is essentially unique. We then call $\k^*$ a
\dfnterm{companion} of $\k$.

The paper \cite{KPT} contains many examples of excellent pairs
$(\k ,\k^*)$ including, e.g., $\k =$ $\a$-free hypergraphs of a given type,
metric spaces with distances in a given countable additive subsemigroup of
$\bbR^+$ (see also \cite{N1} here), vector spaces over a given finite field,
etc., with corresponding companions
$\k^*=$ ordered graphs, ordered $\a$-free hypergraphs of a given type,
ordered metric spaces with distances in a given countable additive
subsemigroup of $\bbR^+$, lexicographically ordered vector spaces over a
given finite field (i.e., with an ordering induced lexicographically by an
arbitrary ordering of a basis), etc. There are also examples of $\k$ that
have no companion (see \cite{J}, \cite{N2}, \cite{KS}), including the class of finite distributive lattices (see \cite{KS}).

If $\k$ is a Fra\"{i}ss\'e class with $\bfk ={\rm Flim}(\k )$ and $\k^*$
is a reasonable order expansion of $\k$, the compact space $X_{\k^*}$ is
0-dimensional, i.e., has a basis consisting of clopen sets.
For each finite substructure $\bfa\subseteq\bfk$ and each $\k^*$-admissible ordering $<$ on $\bfa$,
let
\[
N_{\bfa ,<}:=\{{\prec}\in X_{\k^*}\st{\prec} \restrict A=\ord\}.
\]
This is a clopen basis in $X_{\k^*}$ and the class of the sets $N_{\bfa ,<}$ generates the Borel $\sigma$-algebra of $X_{\k^*}$.

Thus, if $\mu$ is a measure on $X_{\k^*}$, then $\mu$ is completely
determined by the values $\mu (N_{\bfa ,<})$ for $\bfa\subseteq\bfk$ and
$\langle \bfa ,< \rangle\in\k^*$.  These satisfy the obvious relations:
\begin{enumerate}
\item[(i)] For $\bfa\in\k$, $\bfa\subseteq\bfk$,
\[
\sum\big\{\mu (N_{\bfa ,<})\st\langle\bfa ,<\rangle\in\k^*\big\}=1\,.
\]
\item[(ii)] For $\bfa ,\bfb\in\k$, $\bfa\subseteq\bfb\subseteq\bfk$,
\[
\mu (N_{\bfa ,<})=\sum\big\{\mu (N_{\bfb ,<'})\st\langle\bfb ,<'
\rangle\in\k^*,\; \ord\subseteq \ord'\big\}.
\]
\end{enumerate}
Moreover any map $N_{\bfa ,<}\mapsto\mu (N_{\bfa ,<})\in [0,1]$ 
that satisfies these conditions extends uniquely to a measure on $X_{\k^*}$.

In fact, if $\bfa_0\subseteq\bfa_1\subseteq\cdots$ is a sequence of
finite substructures of $\bfk$ with $\bfk =\bigcup_n\bfa_n$, the sets 
$N_{\bfa_n,<}$, for $\langle \bfa_n,<\rangle\in\k^*$, form in an
obvious way a finite branching tree (where the children of $N_{\bfa_n
,<}$ are all $N_{\bfa_{n+1},<'}$ where $\ord\subseteq \ord'$) and a measure
$\mu$ on $X_{\k^*}$ is completely determined by the values $\mu
(N_{\bfa_n,<})$.  Conversely, any map $N_{\bfa_n,<}\mapsto\mu
(N_{\bfa_n,<})\in [0,1]$ such that 
\[
\sum\big\{\mu (N_{\bfa_0,<})\st
\langle\bfa_0,<\rangle\in\k^*\big\}=1
\] 
and 
\[
\mu (N_{\bfa_n,<})=\sum\big\{\mu (N_{\bfa_{n+1},<'})\st\langle\bfa_{n+1},
<'\rangle\in\k^*,\;\ord\subseteq
\ord'\big\}
\]
extends uniquely to a measure on $X_{\k^*}$.

Let $\k$ be a Fra\"{i}ss\'e class in $L$ and $\k^*$ a reasonable
order expansion of $\k$ in $L^*$.  A \dfnterm{consistent random} 
$\k^*$-\dfnterm{admissible ordering} on $\k$ is a map that assigns
to each $\bfa\in\k$ a
probability measure $\mu_\bfa$ on the set of $\k^*$-admissible
orderings on $\bfa$ that is isomorphism invariant (i.e., if $\pi
\colon\bfa\to\bfb$ is an isomorphism, then $\pi_*\mu_\bfa =\mu_\bfb )$
and satisfies the following for each $\bfa\subseteq\bfb$ in $\k$ and
$\Seq{\bfa ,<}\in\k^*$:
\[
\mu_\bfa (<)=\sum\big\{\mu_\bfb (<')\st\langle\bfb ,<'\rangle\in\k^*,\;
\ord\subseteq \ord'\big\},
\]
where we put $\mu_\bfa (<):=\mu_\bfa (\{<\})$.

It is clear that if $(\mu_\bfa )_{\bfa\in\k}$ is a consistent random
$\k^*$-admissible ordering on $\k$, then we can define a $\gp$-invariant
measure $\mu$ on $X_{\k^*}$, where $\gp:={\rm Aut}(\bfk )$, $\bfk :=
{\rm Flim}(\k )$, as follows:
\[
\mu (N_{\bfa ,<}):=\mu_\bfa (<).
\]
Conversely, given such a $\gp$-invariant measure $\mu$ on $X_{\k^*}$,
we can define a consistent random $\k^*$-admissible ordering of $\k$
by
\[
\mu_\bfa (<):=\mu (N_{\bfa',<'}),
\]
where $\bfa'\subseteq\bfk$ and $\langle\bfa ,<\rangle\cong\langle
\bfa' ,<'\rangle$.
By the ultrahomogeneity of $\bfk$, this is
well-defined.

Thus $\gp$-invariant measures on $X_{\k^*}$ can be identified with
consistent random $\k^*$-admissible orderings on $\k$.

We then have, using \rref p.Char/:

\procl p.amen
Let $(\k ,\k^*)$ be an excellent pair.  Let $\bfk :={\rm Flim}(\k
)$.  Then $\gp:={\rm Aut}(\bfk )$ is amenable iff $\k$ admits a consistent random $\k^*$-admissible ordering. Moreover, $\gp$ is uniquely ergodic iff $\k$ admits a unique consistent random $\k^*$-admissible ordering.
\endprocl

If $(\k ,\k^*)$ is an excellent pair for which {\it every} linear
ordering on each given $\bfa\in\k$ is $\k^*$-admissible (in this case,
we write $\k^*=\k *{\mathcal{LO}}$), then there is an obvious Aut$(\bfk
)$-invariant measure $\mu$ on $X_{\k^*}$, which in this case is the space
of {\it all} linear orderings on $K$, given by
\[
\mu (N_{\bfa ,<}):=\frac{1}{n!}\,,
\]
where $n :=|A|$; we call this the \dfnterm{uniform measure} on $X_{\k^*}$.
Thus we have the following:
\begin{pro}
Let $(\k ,\k^* )$ be an excellent pair such that $\k^*=\k *{\mathcal{LO}}$.
Then ${\rm Aut}(\bfk )$ is amenable for $\bfk :={\rm Flim}(\k )$. 
\end{pro}

As we mentioned earlier in this section, examples of classes $\k$ of finite
structures for which $(\k ,\k *{\mathcal{LO}})$ is an excellent pair
include the following (see \cite[Sections 6 and 8]{KPT}): pure sets,
$\a$-free hypergraphs of a given type, metric spaces with (non-zero) distances in a
given countable additive subsemigroup of $\bbR^+$, etc.

\bsection{Order Transitivity and Unique Ergodicity}{s.ot}

We shall discuss here a simple criterion for unique ergodicity and use it to
provide our first examples of uniquely ergodic automorphism groups.
Given an excellent pair $(\k ,\k^*)$, a structure $\bfa\in\k$ is 
$\k^*$-\dfnterm{order transitive} if Aut$(\bfa )$ acts transitively
on the set of $\k^*$-admissible orderings on $\bfa$.  We now have:
\begin{pro}
Let $(\k ,\k^*)$ be an excellent pair, let $\bfk := {\rm Flim}(\k)$ and
assume that $\gp:= {\rm Aut}(\bfk)$ is amenable.
If the class of order-transitive structures in $\k$ is 
cofinal, i.e., for every $\bfa\in\k$, there is an order-transitive
$\bfb\in\k$ with $\bfa\subseteq\bfb$, then $\gp$ is uniquely ergodic.
\end{pro}

\proof  By \rref t.KPT/ and \rref p.Char/, it is enough to show that $X_{\k^*}$ is
uniquely ergodic.  Since $\gp={\rm Aut}(\bfk )$ is amenable, there is
a $\gp$-invariant measure on $X_{\k^*}$.  Fix any such measure $\mu$.
By hypothesis, there is a sequence $\bfa_0\subseteq\bfa_1\subseteq
\cdots$ of finite substructures of $\bfk$ with $\bigcup_n\bfa_n=\bfk$
and each $\bfa_n$ order transitive.  It follows that for each
$\k^*$-admissible order $<$ on $\bfa_n$, we have $\mu (N_{\bfa_n,<})=
{1}/{k_n}$, where $k_n$ is the cardinality of the set of
$\k^*$-admissible orders on $\bfa_n$.  Thus $\mu$ is uniquely
determined.
\Qed
\medskip

Here are some examples of Hrushovski classes of finite structures satisfying the
hypothesis of the previous proposition:
\begin{enumerate}
\item[i)] $\k =$ pure sets, $\k^*=$ linear orderings.  Then $\bfk
={\rm Flim}(\k )=\bbN$ and Aut$(\bfk )=
S_\infty$. So $S_\infty$ is uniquely ergodic (Glasner and Weiss \cite{GW1}).
\item[ii)] $\k =$ equivalence relations, $\k^*=$ equivalence relations
with convex orderings, 
i.e., orderings in which each equivalence class is convex (see \cite[Section
6]{KPT}).
Then $\bfk =$ the equivalence relation on $\bbN$ that has infinitely
many classes, each infinite, so Aut$(\bfk )\cong S_\infty\ltimes
S^\bbN_\infty$ (where $S_\infty$ acts on $S^\bbN_\infty$ by shift)
is uniquely ergodic.
\item[iii)]  Let $\bfk :=\langle\bbN^{<\bbN},\emptyset ,p\rangle$, where
$\emptyset$ is the empty sequence and $p$ is the prefix map,
$p(s):=s\restrict (n-1)$ for $s\in\bbN^n$ for $n>0$ and
$p(\emptyset ):=\emptyset$. 
Then $\bfk$ is a Fra\"{i}ss\'e structure.  Let
$\k :={\rm Age}(\bfk )$.  The elements of $\k$ are the 
structures isomorphic to finite subtrees $T$ of $\bbN^{<\bbN}$, i.e.,
finite subsets $T\subseteq\bbN^{<\bbN}$ containing $\emptyset$ and having
the property that if $s\in\bbN^n$, $n>0$, and $s\in T$, then for every $m<n$, we
have
$s\restrict m\in T$.  We define the class $\k^*$ by saying
that $<$ is a $\k^*$-admissible ordering on $T$ as above if $s,t\in
T$, $s\restrict m=t\restrict m$ (for $m<{\rm length}(s),{\rm length}(t)$) and $s<u<t$,
imply $u\restrict m=s\restrict m$.  Then $(\k ,\k^*)$ is excellent and $\k$ is a
Hrushovski class (see \cite[Section 6.11]{KR} and \cite[Chapter 2]{N1}).  
Moreover, if $T_n:=\{s\in
\bbN^{<\bbN}\st{\rm length}(s)\leq n$ and $\forall i<{\rm length}(s)\;
(s_i<n)\}$, then $\{T_n \st n \ge 0\}$
is a cofinal class of order-transitive structures in $\k$.  Thus 
${\rm Aut}(\bfk )\cong{\rm Aut}(T_\infty )$, where $T_\infty$ is
the rooted $\aleph_0$-regular tree, is uniquely ergodic.
\item[iv)] More generally, let $S\subseteq (0,\infty )$ be countable,
let $\calu_S:=$ the class of finite ultrametric spaces with distances
in $S$ and $\calu^*_S:=$ the class of all convexly ordered finite
ultrametric spaces, where an ordering is \dfnterm{convex} if metric balls are
convex.  Then $(\calu_S,\calu^*_S)$ is excellent and $\calu_S$ is
Hrushovski. Moreover the order-transitive $\bfa\in\calu_S$ are cofinal
(see \cite[Chapter 2] {N1}).  Here $\bfu^{ult}_S:={\rm Flim}(\calu_S)$ is the
\dfnterm{Urysohn
ultrametric space} with distances in $S$ and thus Aut$(\bfu^{ult}_S)=
{\rm Iso}(\bfu^{ult}_S)$ is uniquely ergodic.  (The case $S=\{2^{-n}\st
n\in\bbN\}$ corresponds to the previous example.)  
\item[v)] Let $F$ be a finite field and let $\k_F$ be the class of
finite-dimensional vector spaces over $F$.  Fix an ordering of $F$ in
which $0\in F$ is least.  Let $\k^*_F$ be the class of naturally
ordered vector spaces, where a \dfnterm{natural order} is one induced
lexicographically by an ordering of a basis.  Then $(\k_F,
\k^*_F)$ is an excellent pair and $\k_F$ is Hrushovski (see \cite[Section
6]{KPT}).
Clearly every $\bfa\in \k_F$ is order transitive.  Now $\bfv_{\infty
,F}:={\rm Flim}(\k_F )$ is the (countably) infinite-dimensional vector
space over $F$ and Aut$(\bfv_{\infty ,F})={\rm GL}(V_{\infty ,F})$ is the
general linear group of $\bfv_{\infty ,F}$.  Thus ${\rm GL}(\bfv_{\infty
,F})$ is uniquely ergodic.
\end{enumerate}

\procl r.extend
Let $\bfv\subseteq\bfv_{\infty,F}$ be a
finite-dimensional vector space over $F$.  Then the number of 
$\k^*_F$-admissible orders on $\bfv$ is equal to $|{\rm GL}(\bfv )|$ and if
$\mu$ is the unique invariant measure, then $\mu (N_{\bfv ,<})=
\tfrac{1}{|GL(\bfv )|}$.  Thus if $\bfv\subseteq\bfw\subseteq
\bfv_{\infty,F}$ and $<$ is a $\k^*_F$-admissible order on $\bfv$,
then the number of $\k^*_F$-admissible orders on $\bfw$
that extend $<$ is equal to $\tfrac{|{\rm GL}(\bfw )|}{|{\rm
GL}(\bfv)|}$.
\endprocl

\bsection{A Quantitative Ordering Property and Unique Ergodicity}{s.qo}

We formulate here a quantitative version of the ordering property and show that it implies unique ergodicity for automorphism groups.

Let $\k$ be a Fra\"{i}ss\'e class and $\k^*$ be an order expansion of $\k$.
We say that $\k^*$ satisfies the \dfnterm{quantitative ordering property
(QOP)} if there is an isomorphism-invariant map that assigns to each
structure $\bfa^* =\langle \bfa ,  <_\bfa \rangle\in \k^*$ a real number
$\rho (\bfa^* )\in [0,1]$ such that for every $\bfa \in \k^*$ and every
$\epsilon > 0$, there is a $\bfb = \bfb(\bfa, \epsilon)\in \k$ and a
nonempty set of embeddings $E = E(\bfa, \epsilon)$ of $\bfa$ into $\bfb$
with the property that for each pair of $\k^*$-admissible orderings  $<_\bfa $ of
$\bfa$ and   $<_\bfb $ of $\bfb$, the proportion of embeddings in $E$ that
preserve $<_\bfa$ and $<_\bfb$ is equal to $\rho (\langle\bfa , <_\bfa\rangle)$,
within $\epsilon$.

There is also a slight variation of this property, which we denote by
$\textrm{QOP}^*$, that reads as follows: Let $\k$ be a Fra\"{i}ss\'e class
and $\k^*$ be an order expansion of $\k$.  We say that $\k^*$ satisfies the
\dfnterm{QOP*} if there is an isomorphism-invariant map that assigns to
each structure $\bfa^* =\langle \bfa ,  <_\bfa \rangle\in \k^*$ a real
number $\rho (\bfa^* )\in [0,1]$ such that for every $\bfa^* =\langle\bfa ,
<_\bfa\rangle \in \k^*$ and every $\epsilon > 0$, there is a $\bfb =
\bfb(\bfa^*, \epsilon)\in \k$ and a nonempty set of embeddings $E =
E(\bfa^*, \epsilon)$ of $\bfa$ into $\bfb$ with the property that for each
$\k^*$-admissible ordering  $<_\bfb $ of $\bfb$, the proportion of
embeddings in $E$ that preserve $<_\bfa$ and $<_\bfb$ is equal to $\rho
(\langle\bfa , <_\bfa\rangle)$, within $\epsilon$.

The QOP implies the $\textrm{QOP}^*$; for Hrushovski classes, they are
equivalent by \rref t.ue/. Note that the QOP does not imply the ordering
property, unless the function $\rho$ above is strictly positive.

We now have the following result, whose proof is related to that of \rref l.TVbd/.

\procl p.QOP
Let $\k$ be a Fra\"{i}ss\'e class and $\k^*$ be a
Fra\"{i}ss\'e class that is a reasonable order expansion of
$\k$. Write $\bfk :={\rm Flim}(\k )$ and $\gp:={\rm Aut}(\bfk )$. If $\gp$ is
amenable and the $\textrm{QOP}^*$ holds for $\k^*$, then the $\gp$-flow $X_{\k^*}$ is uniquely ergodic. If moreover $(\k,\k^*)$ is an excellent pair, then $\gp$ is uniquely ergodic.
\endprocl

\proof  Since $\gp$ is amenable, the $\gp$-flow
$X_{\k^*}$  has an invariant measure $\mu$. Let $(\mu_\bfa)_{\bfa\in\k}$ be
the associated consistent random $\k^*$-ordering. For each $\langle\bfa
,<_\bfa\rangle\in\k^*$, we shall show that
$\mu_\bfa (<_\bfa)=\rho (\langle \bfa ,<_\bfa\rangle)$, where $\rho$ comes
from the $\textrm{QOP}^*$; this shows the uniqueness
of $\mu$.

Fix such $\Seq{\bfa ,<_\bfa}$ and $\epsilon >0$.  Let $\bfb$ and $E$ be
as in the definition of the $\textrm{QOP}^*$.
For each $f\in E$, we have $$
\mu_\bfa (<_\bfa) = \mu_{f(\bfa)}\big(f_*(<_\bfa )\big) = \sum \big\{
\mu_\bfb(<_\bfb)\st\langle\bfb ,\;<_\bfb \rangle
\in\k^*, f_* (<_\bfa)\subseteq \ord_\bfb \big\},
$$so
\begin{align*}
|E|\cdot \mu_\bfa (<_\bfa) &=\sum_{f\in E
}\sum\big\{\mu_\bfb(<_\bfb)\st\langle
\bfb ,<_\bfb\rangle\in\k^*,\;f _*(<_\bfa)\subseteq \ord_\bfb\big\}\\
&=\sum_{\langle \bfb ,<_\bfb\rangle\in\k^*}\sum\big\{\mu_\bfb (<_\bfb)\st f
\in E,\; f_*(<_\bfa)\subseteq \ord_\bfb\big\}\\
&=\sum_{\langle \bfb ,<_\bfb\rangle\in\k^*}\mu_\bfb (<_
\bfb)\cdot \big|\big\{f\in
E\st\ f_* (<_\bfa)\subseteq \ord_\bfb\big\}\big|\,,
\end{align*}
and thus
\[
\mu_\bfa(<_\bfa)=\sum_{\langle B,<_\bfb\rangle\in\k^*}\frac{\big|\big\{f\in
E\st\ f_* (<_\bfa)\subseteq \ord_\bfb\big\}\big|}{|E|}
\cdot\mu_\bfb (<_\bfb)\,.
\]

Since $\sum_{\langle \bfb ,<_\bfb\rangle\in\k^*}\mu_\bfb (<_\bfb)=1$, 
this shows that $\big|\mu_\bfa (<_\bfa)-\rho \big(\langle \bfa
,<_\bfa\rangle\big)\big|<\epsilon$, and
the proof is complete.
\Qed

\medskip
In Sections \ref{s.quant}--\ref{s.metric},
we have seen that many excellent pairs $(\k, \k^*)$ satisfy the QOP and
therefore the uniqueness of consistent random $\k^*$-admissible orderings.
As a sample, we have the following result.

\procl t.ua
The automorphism groups of the random $\a$-free hypergraph of a given type
and the Urysohn space $\bfu_S$ are uniquely ergodic, but not compact nor
extremely amenable.
\endprocl

We remark that the Urysohn space $\bfu$ without any restriction on
distances, which is not a \fra\ structure since it is uncountable,
has an extremely amenable isometry group: see \rref b.Pe2/.

\bsection{Hrushovski Structures}{s.Hru}

Let $\bfk$ be a Hrushovski structure.  Then there is a sequence
$\Delta_0\subseteq \Delta_1\subseteq\cdots$ of
compact subgroups of $\gp:={\rm Aut}(\bfk )$ with $\overline{\bigcup_n
\Delta_n}=\gp$.  We shall now prove a stronger version of this fact that
will be used in the next section.
\medskip

\begin{pro}
Let $\bfk$ be a Hrushovski structure.  Then we can find a sequence
$\bfa_0\subseteq\bfa_1\subseteq\cdots$ of finite substructures of
$\bfk$ with $\bfk =\bigcup_n\bfa_n$ and a sequence of compact
subgroups $\Delta_0\subseteq \Delta_1\subseteq\cdots$ of\/ $\gp:={\rm Aut}(\bfk )$
with $\overline{\bigcup_n\Delta_n}=\gp$ such that for each $n$, $A_n$
is invariant under $\Delta_n$ and $g\in \Delta_n\mapsto
g\restrict A_n\in{\rm Aut}(\bfa_n )$ is a surjection from $\Delta_n$ onto 
{\rm Aut}$(\bfa_n)$.
\end{pro}

\proof  Let $\Lambda_0\subseteq \Lambda_1\subseteq\cdots$ be compact
subgroups of $\gp$ with $\overline{\bigcup_n\Lambda_n}=\gp$.  Fix an enumeration
$\{a_0,a_1,\dots\}=K$.  We shall
construct recursively $\bfa_n$, $\Delta_n$ as above such that for each $n$,
$\Lambda_n\subseteq \Delta_n\subseteq \Lambda_{n'}$ for some $n'\geq n$, and $A_n\supseteq
\{a_0,\dots ,a_n\}$.

We take $\bfa_0$ to be a finite substructure of $\bfk$ that contains
$a_0$ and is closed under $\Lambda_0$.  Such exists by the compactness
of $\Lambda_0$.  Let Aut$(\bfa_0)=\{\varphi_1,\dots ,\varphi_k\}$.  Since
$\bigcup_n\Lambda_n$ is dense in $\gp$, there are $f_1,\dots ,f_k\in \Lambda_M$ for some
large $M\geq 1$ such that $f_i\restrict A_0=\varphi_i$ for $1\leq i\leq k$.  Then
$\Lambda_0\cup\{f_1,\dots ,f_k\}\subseteq \Lambda_M$.  Put
\[
\Delta_0:=\overline{\langle \Lambda_0\cup\{f_1,\dots ,f_k\}\rangle}\,,
\]
which is a compact subgroup of $\Lambda_M$ such that $\Lambda_{0}\subseteq
\Delta_0\subseteq \Lambda_{0'}$, where $0':=M$.  Clearly $A_0$ is invariant under
$\Delta_0$ and the restriction map from $\Delta_0$ to Aut$(\bfa_0)$ is
surjective.

Assume now $\bfa_n$, $\Delta_n$ have been constructed.  To define $\bfa_{n+1}
$, $\Delta_{n+1}$, we proceed as before.  Let $\bfa_{n+1}$ be a finite
substructure of $\bfk$ with $\bfa_n\subseteq\bfa_{n+1}$ such that $a_{n+1}
\in A_{n+1}$ and $A_{n+1}$ is invariant under $\Lambda_m$, where $m:=
{\rm max}\{n+1,n'\}$, so that also $\Lambda_m\supseteq \Delta_n\cup \Lambda_{n+1}$.
Let Aut$(\bfa_{n+1})=\{\psi_1,\dots ,\psi_\ell\}$.  As before, there is
$N\geq m$ and $f_1,\dots ,f_\ell\in \Lambda_N$ such that $f_i\restrict A_{n+1}
=\psi_i$ for $1\leq i\leq\ell$.  Put
\[
\Delta_{n+1}:=\overline{\langle  \Lambda_m\cap\{f_1,\dots ,f_\ell\}\rangle}
\subseteq \Lambda_N,
\]
so that $\Delta_{n+1}$ is compact, $\Delta_{n+1}\supseteq \Delta_n\cup
\Lambda_{n+1}$,
and $\Delta_{n+1}\subseteq \Lambda_{(n+1)'}$, where $(n+1)'=N$.  Finally,
$A_{n+1}$ is invariant under $\Delta_{n+1}$ and the restriction map from
$\Delta_{n+1}$ to Aut$(\bfa_{n+1})$ is surjective.

\Qed

\medskip
Let $\bfk$ be a Hrushovski structure.  A sequence $\Seq{(\bfa_n,\Delta_n
) \st n\ge 0}$ as in the previous theorem will be called
\dfnterm{characteristic}.

\bsection{Equivalence of Unique Ergodicity and the QOP for Hrushovski Structures}{s.hru1}

We shall now consider unique ergodicity in the context of Hrushovski
classes. Note first that if $\k$ is a Hrushovski  class and $\bfk = {\rm
Flim}(\k)$, there is a sequence of finite substructures $\bfa_0 \subseteq
\bfa_1 \subseteq \cdots \subseteq\bfk$ such that $\bigcup_n \bfa_n = \bfk$
and every isomorphism between substructures of $\bfa_n$ extends to an
automorphism of $\bfa_{n+1}$. We now have for every such sequence
$\Seq{\bfa_n}$:

\begin{pro} Assume that the Hrushovski class $\k$ admits a companion $\k^*$. Let $\bfk := {\rm Flim} (\k)$. For
each finite $\bfa \subseteq \bfk$ and each
$\k^*$-admissible ordering $<$
on $\bfa$, let $$\mu_n(<) := \frac{|\{<_n \st \langle \bfa_n , <_n
\rangle \in \k^* ,\; \ord\subseteq \ord_n\}|}{|\{\ord_n\st\langle \bfa_n,
<_n\rangle\in \k^*\}|}$$for any $n$ such that $\bfa\subseteq \bfa_n$. Let
$\calu$ be a non-principal ultrafilter on $\bbN$ and put $$\mu(N_{\bfa ,
<}) := \mu (<) := \lim_{n\to \calu} \mu_n (<).$$Then $\mu$ is a
$\gp$-invariant measure on $X_{\k^*}$, where $\gp := {\rm Aut}(\bfk)$. 
\end{pro}

\proof It is easy to check that $\mu$ defines a measure. We next check its
$\gp$-invariance. Let $g\in \gp$ and let $g(\langle \bfa , <\rangle ) =:
\langle \bfb , <'\rangle$ in order to check that $\mu (< ) = \mu (<')$. Let
$n_0$ be large enough so that $\bfa, \bfb \subseteq \bfa_{n_0}$. Let $n >
n_0$. Then $g$ restricted to $\bfa$ is an isomorphism between substructures
of $\bfa_{n_0}$, so $g$ extends to an automorphism of $\bfa_n$. Clearly
$\mu_n (<) = \mu_n (<')$, so $\mu (<) = \mu (<')$.\hfill$\dashv$

\medskip
This gives us the following formula for every such sequence $\Seq{\bfa_n}$
in the case of unique ergodicity.

\begin{thm} Let $\k$ be a Hrushovski class that admits a companion $\k^*$.
Let $\bfk := {\rm Flim}(\k)$ and $\gp := {\rm Aut}(\bfk)$. If\/ $\gp$ is
uniquely ergodic, then the unique $\gp$-invariant measure $\mu$ on
$X_{\k^*}$ is given by $$\mu(N_{\bfa , <}) := \lim_{n\to\infty} \mu_n(<).$$
\end{thm}

\proof This follows from the preceding proposition and the fact that for
any bounded sequence $\Seq{a_n}$ of reals, $\lim_n a_n = a$ iff for
every non-principal ultrafilter $\calu$ on $\bbN$, we have
$\lim_{n\to\calu} a_n = a$.\Qed

\medskip
We shall next see that for Hrushovski classes $\k$ admitting a companion
$\k^*$, unique ergodicity for $\gp:={\rm Aut}(\bfk )$, with
$\bfk :={\rm Flim}(\k )$, is actually equivalent 
to the quantitative ordering property for $\k^*$.

\procl t.ue
Let $\k$ be a Hrushovski class, $\k^*$ be a
Fra\"{i}ss\'e class that is a reasonable order expansion of
$\k$, and let $\bfk :={\rm Flim}(\k )$ and $\gp:={\rm Aut}(\bfk )$. Then the following are equivalent:
\begin{enumerate}
\item[(i)]  The $\gp$-flow $X_{\k^*}$ is uniquely ergodic.
\item[(ii)] There is an isomorphism-invariant map $\rho : \k^* \to [0, 1]$ 
such that for every $\bfa\in \k$ and every $\epsilon >0$, there is
$\bfb\in\k$ with $\bfb \supseteq\bfa$ such that for every $\langle\bfa
,<_\bfa\rangle\in \k^*$ and $\langle\bfb ,<_\bfb\rangle\in\k^*$, the
proportion of automorphisms $\pi$ of $\bfb$ such that $\pi_*(<_\bfa)
\subseteq {<_\bfb}$ is
equal to $\rho (\langle\bfa , <_\bfa\rangle)$,  within $\epsilon$.
\item[(iii)] There is an isomorphism-invariant map $\rho : \k^* \to [0, 1]$
such that for every $\bfa^* =\langle\bfa ,<_\bfa\rangle\in\k^*$ and every
$\epsilon >0$, there is $\bfb \in \k$ with $\bfb\supseteq\bfa$ such that for
every $\langle \bfb ,<_\bfb \rangle\in\k$, the proportion of automorphisms
$\pi$ of $\bfb$ such that $\pi_*(<_\bfa) \subseteq {<_\bfb}$ is
equal to $\rho (\langle\bfa , <_\bfa\rangle)$,  within $\epsilon$.
\item[(iv)] (QOP) There is an isomorphism-invariant map $\rho : \k^* \to [0, 1]$
such that for every $\bfa\in \k$ and each $\epsilon > 0$, there is a $\bfb
\in \k$ and a nonempty set of embeddings $E$ of $\bfa$ into
$\bfb$ with the property that for each $\k^*$-admissible ordering $<_\bfa $
of $\bfa$ and each $\k^*$-admissible ordering $<_\bfb $ of $\bfb$, the
proportion of embeddings in $E$ that preserve $<_\bfa, <_\bfb$
is equal to $\rho (\langle\bfa , <_\bfa\rangle)$, within $\epsilon$.
\item[(v)] $(QOP^*)$ There is an isomorphism-invariant map $\rho : \k^*
\to [0, 1]$ such that for every $\bfa^* = \langle\bfa , <_\bfa\rangle\in
\k^*$ and each $\epsilon > 0$, there is a $\bfb \in \k$ and a nonempty set
of embeddings $E$ of $\bfa$ into $\bfb$ with the property that
for each $\k^*$-admissible ordering $<_\bfb $ of $\bfb$, the proportion of
embeddings in $E$ that preserve $<_\bfa, <_\bfb$ is equal to
$\rho (\langle\bfa , <_\bfa\rangle)$, within $\epsilon$.
\end{enumerate}
Moreover, if $\k^*$ has the {\rm OP}, then (i)--(v) are equivalent to
\begin{enumerate}
\item[(vi)] The same as (ii), but with $\rho$ strictly positive.
\item[(vii)] The same as (iii), but with $\rho$ strictly positive.
\item[(viii)] The same as (iv), but with $\rho$ strictly positive.
\item[(ix)] The same as (v), but with $\rho$ strictly positive.
\end{enumerate}
Finally, if $(\k ,\k^*)$ is an excellent pair, then (i)--(ix) are
equivalent to
\begin{enumerate}
\item[(x)] $\gp$ is uniquely ergodic.
\end{enumerate}
\endprocl

\proof  (i) $\Rightarrow$ (ii).  Let $\mu$ be the unique
invariant measure for the $\gp$-flow $X_{\k^*}$.  Put
\[
\rho (\langle \bfa ,<\rangle):=\mu (N_{\bfa ,<})
\]
for any $\bfa\in\k$, $\bfa\subseteq\bfk$ and $\langle\bfa ,<\rangle\in
\k^*$.  This extends in an obvious way to an isomorphism-invariant
map on all of $\k^*$, also denoted by $\rho$ (since each $\bfa\in\k$ has an
isomorphic copy contained in $\bfk$ and any two such copies are
isomorphic via an automorphism of $\bfk$).

Consider now a characteristic sequence $\Seq{(\bfa_n,\Delta_n) \st n \ge
0}$ as in \rref s.Hru/. Write $\mu_n$ for the Haar measure of $\Delta_n$.
Then
\[
A_n(f)(x):=\int_{\Delta_n}f(g^{-1}\cdot x)\,d\mu_n(g)
\]
converges uniformly to $\int f\,d\mu$ for every $f\in C(X_{\k^*})$ by
\rref t.aver/.  For $\bfa^*:={\langle\bfa ,<\rangle}\in\k^*$ with $\bfa\subseteq
\bfk$, denote by $\bfone_{\bfa^*}$ the indicator function of the set $
\{{\prec}\in X_{\k^*}\st \ord\subseteq{\prec}\}$.
Then  $\bfone_{\bfa^*}$ is continuous, so
\[
\int_{\Delta_n}\bfone_{\bfa^*}(g^{-1}\cdot x)\,d\mu_n(g)\to\int
\bfone_{\bfa^*}\,d\mu
=\rho (\bfa ,<)
\]
uniformly.  For $x={\prec} \in X_{\k^*}$, the left-hand side is
\[
\int_{\Delta_n}\bfone_{\bfa^*}(g^{-1}\cdot x)\,d\mu_n(g)
=\mu_n\big(\{g\in  \Delta_n\st g_*(<)\subseteq {\prec}\}\big).
\]

Find $n_{\bfa^*}$ large enough so that
\rlabel e.1
{\big|\mu_n\big(\{g\in \Delta_n\st g_*(<)\subseteq{\prec}\}\big)-\rho
(\bfa^* )\big|
<\epsilon 
}
for {\it all} ${\prec}\in X_{\k^*}$ and $n\geq n_{\bfa^*}$.  Since there are
only finitely many $\langle\bfa ,<\rangle\in\k^*$ (where $\bfa$ is
fixed), we can find $N$ large enough
so that $\bfa\subseteq\bfa_N$ and \rref e.1/ holds for {\it
all} $<$ with $\bfa^*=\langle\bfa, <\rangle\in\k^*$, ${\prec}\in X_{\k^*}$, and
$n\geq N$.  Take then $\bfb :=\bfa_N$.  Let $\langle\bfa ,<\rangle
\in\k^*$, let $<'$ be $\k^*$-admissible for $\bfb$, and let ${\prec}
\in X_{\k^*}$ extend $<'$.  We only have to check that
\[
\mu_N\big(\{g\in \Delta_N\st g_*(<)\subseteq{\prec}\}\big)
= \frac{|\{\pi\in{\rm Aut}(\bfb )\st\pi_* (<)\subseteq{\prec}\}|}
{|{\rm Aut}(\bfb )|}.
\]
Indeed, let $\varphi\colon \Delta_N\to{\rm Aut}(\bfb )$ be the epimorphism
$\varphi (g):=g\restrict B$.  Then
\begin{align*}
\mu_N\big(\{g\in \Delta_N\st g_*(<)\subseteq{\prec}\}\big)
=&\ \mu_N\big(\{g\in \Delta_N\st\varphi (g)_*(<)\subseteq{\prec}\}\big)\\
=&\ \mu_N\big({\rm ker}(\varphi )\big)\cdot |\{\pi\in{\rm Aut}(\bfb )\st
\pi_* (<)\subseteq{\prec}\}|.
\end{align*}
But clearly $\mu_N\big({\rm ker}(\varphi )\big)=\tfrac{1}{|{\rm Aut}(\bfb
)|}$, so we are done.

\smallskip
(ii) $\Rightarrow$ (iii) is obvious.

\smallskip
(iii) $\Rightarrow$ (i). This is similar to the proof of \rref p.QOP/. Since $\gp$ is amenable, the $\gp$-flow
$X_{\k^*}$  has an invariant measure.  Let $\mu$ be any such measure.
For each $\Seq{\bfa ,<_\bfa}\in\k^*$ with $\bfa\subseteq\bfk$, we shall
show that $\mu (N_{\bfa ,<_\bfa})=\rho (\langle \bfa ,<_\bfa\rangle)$,
which shows the uniqueness of $\mu$.

Fix such $\Seq{\bfa ,<_\bfa}$ and $\epsilon >0$.  Let $\bfb$ be as in (iii).
For each $\pi\in{\rm Aut}(\bfb )$, we have
\[
\mu (N_{\bfa ,<_\bfa})=\sum\big\{\mu(N_{\bfb ,<_\bfb})\st\langle\bfb ,<_\bfb\rangle
\in\k^*,\;\pi_* (<_\bfa)\subseteq \ord_\bfb\big\}
\]
since $\mu (N_{\bfa ,<_\bfa})=\mu (N_{\pi_* (\bfa ),\pi_* (<_\bfa)})$ and $N_{\pi_*
(\bfa ),\pi_* (<_\bfa)}$ is the disjoint union of the sets $N_{\bfb
,<_\bfb}$ with $\langle \bfb ,<_\bfb\rangle\in\k^*$ and $\pi_*
(<_\bfa)\subseteq \ord_\bfb$.  So
\begin{align*}
|{\rm Aut}(\bfb )|\cdot \mu (N_{\bfa ,<_\bfa}) &=\sum_{\pi\in{\rm Aut}
(\bfb )}\mu (N_{\bfa ,<_\bfa})\\
&=\sum_{\pi\in{\rm Aut}(\bfb )}\sum\big\{\mu (N_{\bfb ,<_\bfb})\st\langle
\bfb ,<_\bfb\rangle\in\k^*,\;\pi_* (<_\bfa)\subseteq \ord_\bfb\big\}\\
&=\sum_{\langle \bfb ,<_\bfb\rangle\in\k^*}\sum\big\{\mu (N_{\bfb ,<_\bfb})\st\pi
\in{\rm Aut}(\bfb ),\;\pi_* (<_\bfa)\subseteq \ord_\bfb\big\}\\
&=\sum_{\langle \bfb ,<_\bfb\rangle\in\k^*}\mu (N_{\bfb ,<_\bfb})\cdot |\{\pi\in
{\rm Aut}(\bfb )\st\pi_* (<_\bfa)\subseteq \ord_\bfb\}|,
\end{align*}
whence
\[
\mu (N_{\bfa ,<_\bfa})=\sum_{\langle \bfb,<_\bfb\rangle\in\k^*}\frac{|\{\pi\in
{\rm Aut}(\bfb )\st\pi_* (<_\bfa)\subseteq \ord_\bfb\}|}{|{\rm Aut}(\bfb )|}
\cdot\mu (N_{\bfb ,<_\bfb}).
\]

Since $\sum_{\langle \bfb ,<_\bfb\rangle\in\k^*}\mu (N_{\bfb ,<_\bfb})=1$, 
this shows that we have $|\mu (N_{\bfa ,<_\bfa})-\rho (\langle \bfa ,<_\bfa\rangle)|<\epsilon$.

\smallskip
(ii) $\Rightarrow$ (iv). Let $\rho$ be as in (ii). Given $\bfa$ and
$\epsilon > 0$, let $\bfb$ again be as in (ii). Let $E$ consist
of the restrictions to $\bfa$ of all the automorphisms of $\bfb$. For $f \in
E$, let $A_f := \{ \pi \in {\rm Aut} (\bfb )\st \pi \restrict
A = f\}$. It is clearly enough to show that $|A_f|$ is independent of $f$.
So take $f_1, f_2 \in E$ and fix $\pi_1 \in A_{f_1}$, $\pi_2
\in A_{f_2}$. It is enough to show that for $\sigma := \pi_2\circ\pi_1^{-1}$,
we have $\sigma\circ A_{f_1} \subseteq A_{f_2}$. Indeed, let $\pi\in
A_{f_1}$. Then for $a\in A$, we have
$\sigma\circ\pi (a) = \sigma\big(f_1
(a)\big) = \sigma\big(\pi_1 (a)\big)= \pi_2 (a) = f_2 (a)$, i.e.,
$\sigma\circ\pi\in A_{f_2}$.

\smallskip
(iv) $\Rightarrow$ (v) is obvious.

\smallskip
(v) $\Rightarrow$ (i). This follows by \rref p.QOP/.

\smallskip
If $\k^*$ satisfies the OP, then $X_{\k^*}$ is a minimal $\gp$-flow and thus
the (closed) support of any invariant measure is equal to $X_{\k^*}$,
whence in all of (ii)--(v), we may take $\rho$ to be strictly positive, i.e., 
(vi)--(ix) hold.

Finally, the equivalence with (x) follows from \rref p.Char/ and \rref t.KPT/.
\Qed

\procl r.strengthening
In \rref t.ue/, consider the following strengthening of (ii):
\begin{enumerate}
\item[(ii)$'$] There is an isomorphism-invariant map $\rho\colon \k^* \to
[0,1]$
such that for every
$\bfa\in\k$ and $\epsilon >0$, there is $\bfb\in\k$ with $\bfb\supseteq\bfa$ such that for
$\bfa_0\subseteq\bfa$, $\langle\bfa_0,<_0\rangle\in\k^*$, and every 
$\k^*$-admissible ordering $<_\bfb$ for $\bfb$, we have
\[
\left |\frac{|\{\pi\in{\rm Aut}(\bfb )\st\pi_* (<_0)\subseteq 
\ord_\bfb\}|}{|{\rm Aut}(\bfb )|}-\rho \big(\langle \bfa_0,<_0\rangle\big)\right |<\epsilon .
\]
\end{enumerate}
\endprocl

It is easy to check that the proof of \rref t.ue/ also shows that for a 
Hrushovski class $\k$, (ii)$'$ is equivalent to (ii).  However, we
can see, without assuming that $\k$ is Hrushovski, that (ii)$'
\Rightarrow$ (i) and thus when $(\k ,\k^*)$ is an excellent pair,
(ii)$'$ implies that $\gp:={\rm Aut}(\bfk )$ is uniquely ergodic.

We shall check that $\mu_\bfa (<):=\rho (\bfa ,<)$ for $\langle\bfa
,<\rangle\in\k^*$ is a consistent random $\k^*$-admissible ordering
on $\k$.  This shows that the $\gp$-flow $X_{\k^*}$ admits an
invariant measure and the argument in \rref t.ue/ (iii) $\Rightarrow$ (i)
shows that it is uniquely ergodic.

First note that $$ \forall \bfa\in\k \quad \sum\big\{\rho (\bfa^*)\st
\bfa^*\restrict L=\bfa\big\}=1,$$since (ii)$'$ implies that for all
$\epsilon$, this sum is equal to 1 within $|A|!\epsilon$.

Next fix $\bfa_0\subseteq\bfa$ in $\k$ and $(\bfa_0,<_0)\in\k^*$ in 
order to show that
\[
\rho \big(\langle \bfa_0,<_0\rangle\big)=\sum\big\{\rho (\bfa ,<)\st\langle\bfa ,<\rangle
\in\k^*,\; \ord_0\subseteq {<}\big\}.
\]
Let $\epsilon >0$ and let $\bfb$ be as in (ii)$'$.

Note that for $\langle\bfb ,<_\bfb\rangle\in\k^*$,
\begin{align*}
\big\{\pi\in{\rm Aut}(\bfb )\st\pi_* (<_0)\subseteq \ord_\bfb\big\}
=&\ \big\{\pi\in{\rm Aut}(\bfb )\st\ord_0\subseteq\pi^*(<_\bfb)\big\}\\
=&\ \bigsqcup_{\substack{\langle\bfa ,<\rangle\in\k^*\\ \ord_0\subseteq
\ord}}
\big\{\pi\in{\rm Aut}(\bfb )\st\pi_* (<)\subseteq {<_\bfb}\big\},
\end{align*}
so
\[
\frac{\big|\big\{\pi\in{\rm Aut}(\bfb )\st\pi_* (<_0)\subseteq
\ord_\bfb\big\}\big|}{|
{\rm Aut}(\bfb )|}=
\]
\begin{gather}
\sum\left\{\frac{\big|\big\{\pi\in{\rm Aut}(\bfb )\st
\pi_* (<)\subseteq \ord_\bfb\big\}\big|}{|{\rm Aut}(\bfb )|}
\st\langle\bfa ,<\rangle\in\k^*,\; \ord_0\subseteq \ord\right\}
\,.
\tag{$*$}
\end{gather}

Let
\[
a:=\sum\big\{\rho (\langle \bfa ,<\rangle)\st\langle\bfa
,<\rangle\in\k^*,\;\ord_0
\subseteq {<}\big\}
\]
and
\[
N:=\big|\big\{{<}\st\langle\bfa ,<\rangle\in\k^*,\;\ord_0\subseteq
{<}\big\}\big|\,.
\]
Then
\[
a-N\epsilon <(*)<a+N\epsilon\,.
\]
Thus
\[
\left |
\frac{\big|\big\{\pi\in{\rm Aut}(\bfb )\st\pi_* (<_0)\subseteq
\ord_\bfb\big\}\big|}{{|\rm Aut}(\bfb)|}
-\sum\big\{\rho (\langle \bfa ,<\rangle)\st\langle\bfa ,<\rangle\in\k^*,\;
\ord_0\subseteq
{<}\big\}\right |<N\epsilon.
\]
But also
\[
\left |\frac{\big|\big\{\pi\in{\rm Aut}(\bfb )\st\pi_* (<_0)\subseteq
\ord_\bfb\big\}\big|}
{|{\rm Aut}(\bfb )|}-\rho (\langle \bfa_0,<_0\rangle)\right |<\epsilon ,
\]
so, taking $\epsilon\to 0$, we see that
\[
\rho (\bfa_0,<_0)=\sum\big\{\rho (\langle \bfa ,<\rangle)\st\langle\bfa ,<\rangle
\in\k^*,\; \ord_0\subseteq {<}\big\}.
\]

A similar remark holds if we replace (iv) by the analogous (iv)$'$.
\bsection{The Support of the Unique Measure}{s.sup}

We shall show here that in certain situations where unique ergodicity holds, the unique measure is supported by a single orbit, which is actually comeager. 

We first discuss the notion of generic point and orbit. If $\gp$ is a
topological group that acts continuously on a topological space $X$, we
say that $x\in X$ is \dfnterm{generic} if its orbit $\gp\cdot x$ is comeager.
In this case, we also say that $\gp \cdot x$ is a \dfnterm{generic orbit}. Clearly there is at most one generic orbit in any Baire space $X$. We first note the following general fact:

\begin{pro}
Let\/ $\gp$ be a Polish group acting continuously on topological spaces
$X$ and $Y$ that are Hausdorff and Baire. Assume that the action of\/ $\gp$
on $X$ is minimal and $\pi : X\rightarrow Y$ is a continuous surjective
$\gp$-map. If $x_0$ is a generic point for $X$, then $\pi(x_0)$ is a
generic point for $Y$.
\end{pro}

\proof We use arguments similar to those in Appendix A of Melleray and Tsankov
\cite {MT}, although we need to exercise extra care because of our more
general context.

Let $y_0 := \pi (x_0)$. First notice that $\gp\cdot y_0$ is dense in $Y$.
We next verify that $\gp\cdot y_0$ has the Baire Property in $Y$. By the
Nikod\'{y}m Theorem (see \cite[29.14]{K}), it is enough to show that
$\gp\cdot y_0$ can be obtained via the Souslin operation $\a$ applied to
closed sets in $Y$, i.e., can be written in the form $\a_s F_s := \big\{
y\in Y \st\exists \alpha \in \n \ \forall n\  (y\in F_{\alpha \restrict n
})\big\}$, where $\n := \bbN^\bbN$ is the Baire space and $(F_s)_{s\in
\bbN^{<\bbN}}$ is a family of closed sets indexed by the set of finite
sequences from $\bbN$. To exhibit such a representation, let $\rho :\n
\rightarrow \gp$ be a continuous surjection and let $f : \n \rightarrow Y$
be defined by $f(\alpha ) := \rho (\alpha ) \cdot y_0$.  Clearly $f$ is
continuous.
Let $\n_s := \{ \alpha \in \n \st s\subseteq \alpha\}$,  for $s\in
\bbN^{<\bbN}$, be the basic open sets in $\n$. Put $F_s :=
\overline{f[\n_s]}$. Then it is easy to see that $\gp\cdot y_0 = \a_s F_s$.

It follows that if $\gp\cdot y_0$ is not meager, it must be comeager (see
\cite[8.46]{K}). That is, it suffices to show that $\gp\cdot y_0$ is not
meager. Now if it were meager, then there would be 
a sequence $\Seq{V_n}$ of dense open sets in $Y$
with $\bigcap_n V_n \cap \gp\cdot y_0 = \emptyset$. This would imply that
$\bigcap_n \pi^{-1} (V_n) \cap \gp\cdot x_0 =\emptyset$, so it is enough to
show that if $V\subseteq Y$ is dense and open in $Y$, then $\pi^{-1} (V)$
is dense in $X$. To show this, let $W\subseteq X$ be nonempty and open. Let
$\gp_0$ be a countable dense subgroup of $\gp$. By the minimality of $X$,
we have $\gp_0 \cdot W = X$, whence $\gp_0 \cdot \pi (W) = Y$. Then $\pi
(W)$ is not meager in $Y$, thus $\pi (W) \cap V \not= \emptyset$, so $W
\cap \pi^{-1} (V) \not= \emptyset$.  \Qed

\medskip

Specializing to the case of $\gp$-flows, we then have the following: 

\begin{cor}
Let\/ $\gp$ be a Polish group. If the universal minimal flow of\/ $\gp$ has
a generic point, then so does every minimal flow of\/ $\gp$.
\end{cor}

Let us say that a Polish group $\gp$ has the \dfnterm{generic point property} if every minimal $\gp$-flow has a generic point. We note now the following:

\begin{pro} 
Let $(\k, \k^*)$ be an excellent pair. Let $\bfk := {\rm Flim} (\k)$ and
$\gp:= {\rm Aut} (\bfk )$. Then $\gp$ has the generic point property.
\end{pro}

\proof By the previous corollary, it is enough to show that $X_{\k^*}$ has a generic point.
Let
$\bfk^*:=\langle\bfk ,<^*\rangle$ be the Fra\"{i}ss\'e limit of
$\k^*$.  
Clearly the $\gp$-orbit of $<^*$ consists of all $<$ in $X_{\k^*}$ such
that $\langle\bfk ,<\rangle\cong\langle\bfk ,<^*\rangle$, i.e., all
$\langle\bfk ,<\rangle$ that are, up to isomorphism, the
Fra\"{i}ss\'e limit of $\k^*$.  These are characterized by the
following two properties (see \cite[7.1.4]{Ho}):
\begin{enumerate}
\item[(i)] Age$(\langle\bfk ,<\rangle )=\k^*$;
\item[(ii)] Given $\bfa^*\subseteq\bfb^*$ in $\k^*$ and an embedding
$\pi\colon\bfa^*\to\langle\bfk ,<\rangle$, there is an embedding
$\rho\colon\bfb^*\to\langle\bfk ,<\rangle$ extending $\pi$.
\end{enumerate}
Since property (i) is true for all $\ord\in X_{\k^*}$, this orbit consists of
all $\ord \in X_{\k^*}$ that satisfy condition (ii), and this is clearly a
$G_{\delta }$ subset of $X_{\k^*}$. It is also dense by \cite[Section 7]{KPT}.
\Qed

\medskip
We now have the following result.

\begin{thm}
Let $(\k ,\k^*)$ be an excellent pair with
$\k^*=\k *{\mathcal{LO}}$.  Let $\bfk :={\rm Flim}(\k )$ and $\gp:={\rm
Aut}(\bfk )$.
Then the uniform measure on $X_{\k^*}$ is supported by the generic orbit.
In particular, if\/ $\gp$ is uniquely ergodic, then the unique measure in
each minimal flow is supported by the generic orbit.
\end{thm}

\proof Since $\k^*=\k *{\mathcal{LO}}$, it is clear that $X_{\k^*}={\rm LO}
(K)$ is the space of
all linear orderings on $K$. Let
$\bfk^*:=\langle\bfk ,<^*\rangle$ be the Fra\"{i}ss\'e limit of
$\k^*$.  Clearly $\ord^*\in {\rm LO}(K)$. We shall show that the uniform
measure $\mu$ is supported by the generic orbit, which is the $\gp$-orbit
of $<^*$ in the $\gp$-space
${\rm LO}(K)$.

As we have seen in the proof of the preceding proposition, the $\gp$-orbit of $<^*$ consists of all $<$ in ${\rm LO}(K)$ that satisfy the following property:

Given $\bfa^*\subseteq\bfb^*$ in $\k^*$ and an embedding
$\pi\colon\bfa^*\to\langle\bfk ,<\rangle$, there is an embedding
$\rho\colon\bfb^*\to\langle\bfk ,<\rangle$ extending $\pi$.

Since $\k^*$ is countable (up to isomorphism), it is enough to show that
for each given $\bfa^*\subseteq \bfb^*$ in the class $\k^*$, where
$\bfa^*=\langle\bfa, <_\bfa\rangle$ and $\bfb^* =\langle\bfb ,<_\bfb\rangle$,
and any embedding $\pi_\bfa\colon\bfa\to\bfk$, if we let $(\pi_\bfa)_*
(\bfa )=:\bfa_0\subseteq \bfk$ and also $(\pi_\bfa)_* (<_\bfa ) =: \ord_0$,
then the set of all $\ord\in {\rm LO}(K)$ with $\ord_0\subseteq \ord$ for which
there is no embedding $\pi_\bfb\colon\bfb \to\bfk$ extending $\pi_\bfa$
with $(\pi_\bfb)_* (<_\bfb )\subseteq \ord$ is $\mu$-null.

Since $\k$ satisfies the \dfnterm{strong amalgamation property} (see
\cite[Section 2 and 5.3]{KPT} for the definition and this result), for each
$n\geq 1$, we can find $\bfb_1,\dots ,\bfb_n\subseteq\bfk$ such that
$\bfa_0\subseteq\bfb_i$ for $1\leq i\leq n$ and $B_i\cap B_j=A_0$ if $i\neq j$,
and isomorphisms $\pi_i\colon\bfb\to\bfb_i$ extending $\pi_\bfa$. Let $<_i$
be the image of $<_\bfb$ by $\pi_i$.  Thus
$\ord_0\subseteq \ord_i$ for $1\leq i\leq n$.  It is thus enough to show that the
$\mu$-measure of the set of $\ord\in {\rm LO}(K)$ with $\ord\supseteq \ord_0$ but
$\ord\not\supseteq \ord_i,\forall 1\leq i\leq n$, tends to 0 as $n\to
\infty$.

When the language of $\k$ is relational, then it is clear that it is enough
to restrict ourselves to pairs $\bfa^*$, $\bfb^*$ as above where $|B| = |A|
+1$. In this case, we can complete the proof as follows:

Let $A:=\{a_1,\dots ,a_k\}$ and $B=\{a_1,\dots ,a_k,b\}$, where $a_1<_\bfa
a_2<_\bfa a_3<_\bfa\cdots <_\bfa a_k$.  Then one of the following holds:
\[
b<_\bfb a_1,\; a_1<_\bfb b<_\bfb a_2,\; a_2<_\bfb b<_\bfb a_3,\;\dots ,
a_{k-1}<_\bfb b<_\bfb a_k, \mbox{ or } a_k<b .
\]
Assume $a_1<_\bfb b<_\bfb a_2$, the argument being similar
in all the other cases. Let $b_i:=\pi_i(b)$ for $1\leq i\leq n$ and $\pi_\bfa
(a_j)=\overline a_j$ for $1\leq j\leq k$.  Then
\begin{align*}
&\big\{\ord\in {\rm LO}(K)\st \ord_0\subseteq \ord\ \&\ \forall 1\leq i\leq n\;
(\ord_i\not\subseteq \ord)\big\}\\
\subseteq &\big\{\ord\in {\rm LO}(K)\st \ord_0\subseteq \ord\ \&\ \forall 1\leq i\leq
n\;(b_i<\overline a_1\vee\overline a_2<b_i)\big\}\\
\subseteq &\bigcup\big\{N_{\bfc ,<}\st \ord\in {\rm LO}(C)\ \&\ \overline a_1
<\overline a_2\ \&\ \forall 1\leq i\leq n\;(b_i<\overline a_1\vee
\overline a_2<b_i)\big\},
\end{align*}
where $\bfc$ is the substructure of $\bfk$ with universe
$C=\{\overline a_1,\overline a_2,b_1,\dots ,b_n\}$.  Since $\mu
(N_{\bfc ,<})=\tfrac{1}{(n+2)!}$ for each $\ord\in {\rm LO}(C)$, it is enough
to show that
\[
\frac{\big|\big\{<\in {\rm LO}(C)\st\overline a_1<\overline a_2\ \&\ \forall
1\leq i\leq n\;(b_i<\overline a_1\vee\overline a_2<b_i\big\}\big|}{(n+2)!}\to 0
\]
as $n\to\infty$.  But simple counting shows that numerator is
equal to $(n+1)!$, so this ratio is equal to $\tfrac{1}{n+2}\to 0$,
and the proof is complete.

In the general case, we have $A = \{a_1, \dots , a_k\}$ and $B = \{a_1, \dots ,
a_k, b_1, \dots , b_m\}$. Then if $b^i_l := \pi_i (b_l )$ $(1\leq i\leq n,
1\leq l \leq m)$ and $\pi_\bfa (a_j) =: \bar{a}_j$, let $C := \{\bar{a}_1,
\dots , \bar{a}_k\}\cup \bigcup_{i=1}^n\{b^i_1, \dots , b^i_m\}$. Then the
$\mu$-measure of $$\big\{ \ord \in {\rm LO} (K)\st \ord_0 \subseteq \ord \; {\rm and} \;
\forall 1\leq i \leq n\; (\ord_i \nsubseteq \ord)\big\}$$is equal to
$$\frac{\big|\big\{ \ord \in
{\rm LO} (C)\st \ord_0 \subseteq \ord \; {\rm and} \; \forall 1\leq i \leq
n\; (\ord_i
\nsubseteq \ord)\big\}\big|}{| {\rm LO } (C)|},$$so it is enough to show that this ratio
goes to 0 as $n$ goes to infinity. 
This is a consequence of the following lemma.

\begin{lem} Let $X = \{x_1, \dots , x_k\}$ and $Y = \{y_1 , \dots , y_m\}$
be two disjoint sets and let $<_*$ be an ordering of $X \sqcup Y$. Consider
the set $ X + n Y := X\sqcup \bigsqcup^n_{i=1}Y^i$, where $Y^i := \{y^i_1,
\dots , y^i_m\}$, consisting of the union of $X$ and $n$ disjoint copies
of\/
$Y$. Let $<_i$ be the copy of $<_*$ on $X\sqcup Y^i$. Then, for the uniform
probability measure on ${\rm LO}(X+nY)$, the probability that a linear
ordering $<$ on $X+nY$ extends ${<_*} \restrict X$ but $\ord \nsupseteq
\ord_i$ for every $1 \leq i \leq n$, tends to 0 as $n \to \infty$.
\end{lem}

\proof
  This was proved first by Padraic Bartlett with another method.
  One way to generate a uniform ordering on $X+nY$ is to assign each
  element $t\in X+nY$ an independent random variable $U_t$,
  uniform in $[0,1]$. The order $<$ is then induced from $[0,1]$.

  The event that $<$ extends ${<_*}\restrict X$ depends only on $\Seq{U_t \st
  t\in X}$,
  and is equivalent to those being in the right order. Conditioned on
  $\Seq{U_t \st
  t\in X}$, the restrictions of $<$ to $X\sqcup Y^i$ are
  independent, and each has non-zero probability $q = q(\Seq{U_t\st t\in
  X})$ of being $<_i$. Thus the conditional probability that $<$
  extends ${<_*}\restrict X$ but not any $<_i$ is $(1-q)^n$.

  Since $q$ is almost surely non-zero, after taking expectation with
  respect to $\Seq{U_t \st t\in X}$, this tends to $0$ as $n\to\infty$
  by the bounded convergence theorem.  

\Qed

\medskip

Many of the examples of uniquely ergodic automorphism groups that we discussed earlier satisfy the conditions of the preceding theorem, so the unique invariant measure in each minimal flow concentrates on the generic orbit. These include the automorphism groups mentioned in \rref t.ua/.

\bsection{Some Open Problems}{s.problems}

The preceding work suggests a number of open problems.

\begin{question}[Unique Ergodicity Problem] Let $\gp$ be an amenable Polish group with
metrizable universal minimal flow. Is $\gp$ uniquely ergodic?
\end{question}

Recall that the universal minimal flow of a countable group is not metrizable, so infinite countable groups, which, as we mentioned in the introduction, are not uniquely ergodic, do not provide counterexamples (and probably the same holds for non-compact, locally compact groups).

One can even consider a more general version of Question 15.1 for
Polish groups that need not be amenable: If $\Gamma$ is a Polish group
with metrizable universal minimal flow, then does every minimal flow of
$\Gamma$ have at most one invariant measure?

\begin{question}[Generic Point Problem] Let $\gp$ be a Polish group with
metrizable universal minimal flow. Does  $\gp$ have the generic point property?
\end{question}

Finally we have the following stronger version of the first problem.

\begin{question}[Unique Ergodicity-Generic Point Problem] Let $\gp$ be an amenable Polish group with
metrizable universal minimal flow. Is
it true that  $\gp$ is uniquely ergodic and has the generic point property
and moreover for every minimal $\gp$-flow, the unique invariant measure is
supported by the generic orbit?
\end{question}

Recent work of Nguyen Van Th\'e and Tsankov may be relevant to these problems for the case when $\gp$ is the automorphism group of a Fra\"{i}ss\'e structure $\bfk$.

As we mentioned earlier in \rref s.S2/, there are examples of
Fra\"{i}ss\'e classes $\k$ that have no companions at all. However, Nguyen
Van Th\'e \cite{N2} developed a more general notion of expansion $\k^*$
for a given Fra\"{i}ss\'e class $\k$ in a language $L$. Such an expansion
is obtained by taking $L^*$ to be a language obtained from $L$ by adding
not merely a single binary relation symbol $<$, but instead a finite or infinite
(countable) family of relation symbols (of various arities) $(R_i)_{i\in
I}$. It is shown in \cite{N2} that the basic theory of \cite{KPT} goes
through in this more general context, provided the class $\k^*$ is
\dfnterm{precompact}, i.e., every structure in $\k$ has only finitely many
expansions in $\k^*$ (this will be automatically true if one adds only
finitely many symbols to $L$ to form $L^*$). In particular, if $\k$ admits a
precompact companion $\k^*$ (i.e., such an expansion that satisfies the
Ramsey Property and the analog of the Ordering Property in this context,
called the \dfnterm{Expansion Property}), then an analogous metrizable space $X_{\k^*}$ is the universal minimal flow of the automorphism group of the Fra\"{i}ss\'e limit of $\k$. 

Nguyen Van Th\'e and Tsankov (private communication, 2012) have now shown
that for a Fra\"{i}ss\'e class $\k$ with Fra\"{i}ss\'e limit $\bfk$ and
$\gp:= {\rm Aut} (\bfk )$, the following are equivalent:

\begin{itemize}
\item $\gp$ has metrizable universal minimal flow with a comeager orbit,
\item $\k$ admits a precompact expansion with the Ramsey Property and the Expansion Property.
\end{itemize}

Such precompact companions have been computed for: (i)  the class of
local orders, (ii) the age of the Fra\"{i}ss\'e directed graph $S(3)$ (see
\cite{N2}) and (iii) the class of boron tree structures (see \cite {J}).
(In all these cases, the language $L^*$ turns out to be finite.) For (i) it
was shown in \cite{KS} that the corresponding automorphism
group is not amenable and the same has been proved for (ii) and (iii)  by Andrew Zucker. 

Finally, let again $(\k, \k^*)$ be an excellent pair, with $\bfk :={\rm
Flim} (\k)$ and $\gp := {\rm Aut}(\bfk )$ amenable. In all the cases that we
have been able to prove unique ergodicity for $\gp$, it turned out that the
unique $\gp$-invariant measure on $X_{\k^*}$ was the uniform measure given
by $\mu (N_{\bfa , <}) := \frac{1}{k(\bfa)}$, where $k(\bfa) :=
\big|\big\{{<} \st
\langle \bfa , <\rangle \in \k^*\big\}\big|$ is the cardinality of the set of $\k^*$-admissible orderings on $\bfa $. One can ask whether this is a general phenomenon.

\begin{question}
Let $(\k, \k^*)$ be an excellent pair with $\bfk := {\rm Flim}(\k)$ and let
$\gp := {\rm Aut}(\bfk)$ be amenable. Is there a (necessarily unique) $\gp$-invariant measure $\mu$ on $X_{\k^*}$ satisfying $\mu (N_{\bfa , <}) = \frac{1}{k(\bfa)}$?
\end{question}

Notice that this is equivalent to asking the following: Let $(\k, \k^*)$ be as in the previous problem. Is it true that for any $\bfa \subseteq \bfb$ in $\k$, every $\k^*$-admissible ordering on $\bfa$ has the same number of extensions to a $\k^*$-admissible ordering on $\bfb$?

\bigskip
Omer Angel

Department of Mathematics

University of British Columbia

Vancouver, British Columbia V6T 1Z2

\url{angel@math.ubc.ca}

\medskip
Alexander S. Kechris

Department of Mathematics

California Institute of Technology

Pasadena, CA 91125

\url{kechris@caltech.edu}

\medskip
Russell Lyons

Department of Mathematics

831 E 3rd St

Indiana University

Bloomington, IN 47405-7106

\url{rdlyons@indiana.edu}

\url{mypage.iu.edu/~rdlyons}


\begin{thebibliography}{MMM}





\bibitem[BK]{BK} H. Becker and A.S. Kechris (1996). {\it The Descriptive
Set Theory of Polish Group Actions}, Cambridge Univ. Press.

\bibitem[BHV]{BHV} B. Bekka, P. de la Harpe and A. Valette (2008). {\it
Kazhdan's Property {\rm (T)}}, Cambridge Univ. Press.

\bibitem[Fr]{Fr} R. \fra\ (1954).
Sur l'extension aux relations de quelques propri\'et\'es des ordres, 
{\it Ann. Sci. Ecole Norm. Sup.} (3) {\bf 71}, 363--388. 

\bibitem[Gl]{Gl} E. Glasner (2003). {\it Ergodic Theory via Joinings}, Amer.
Math. Soc.

\bibitem[GW1]{GW1} E. Glasner and B. Weiss (2002). Minimal actions of the
group $S(\bbZ )$ of permutations of the integers, {\it Geom. Funct.
Anal.} {\bf 12}, 964--988.

\bibitem[GW2]{GW2} E. Glasner and B. Weiss (2003). The universal minimal
system for the group of homeomorphisms of the Cantor set, {\it
Fund. Math.} {\bf 176}, 277--289.

\bibitem[He1]{He1}B. Herwig (1995). Extending partial isomorphisms on finite
structures, {\it Combinatorica} {\bf 15} (3), 365--371.

\bibitem[He2]{He2} B. Herwig (1998). Extending partial isomorphisms for the
small index property of many $\omega$-categorical structures, {\it Israel
J. Math.} {\bf 107}, 93--123.

\bibitem[Ho] {Ho} W. Hodges (1993). {\it Model Theory}, Cambridge Univ. Press.

\bibitem[Hr]{Hr} E. Hrushovski (1992). Extending partial isomorphisms of
graphs, {\it Combinatorica} {\bf 12}, 411--416.

\bibitem[J]{J} J. Jasi\'nski (2011). Hrushovski and Ramsey Properties of
Classes of Finite Inner Product Structures, Finite Euclidean Metric Spaces,
and Boron Trees, {\it Ph.D. Thesis}, Univ. of Toronto.

\bibitem[K]{K} A.S. Kechris (1995). {\it Classical Descriptive Set
Theory}, Springer.

\bibitem[KPT]{KPT} A.S. Kechris, V.G. Pestov, and S. Todorcevic (2005). 
Fra\"{i}ss\'e limits, Ramsey theory, and topological dynamics of
automorphism groups, {\it Geom. Funct. Anal.} {\bf 15},
106--189.

\bibitem [KR]{KR} A.S. Kechris and C. Rosendal (2007). Turbulence,
amalgamation, and generic automorphisms of homogeneous structures,
{\it Proc. London Math. Soc.} {\bf 94} (3), 302--350.

\bibitem[KS]{KS} A.S. Kechris and M. Soki\'c (2011). Dynamical properties
of the automorphism groups of the random poset and random distributive
lattice, {\it Fund. Math.}, to appear.

\bibitem[McD]{McD}
C. McDiarmid (1989).
\newblock On the method of bounded differences.
\newblock In Siemons, J., editor, {\em Surveys in Combinatorics, 1989}, volume
  141 of {\em London Math. Soc. Lecture Note Ser.}, pages 148--188. Cambridge
  Univ. Press, Cambridge.
\newblock Papers from the Twelfth British Combinatorial Conference held at the
  University of East Anglia, Norwich, 1989.

\bibitem[Mec]{Mec}
M.W. Meckes (2010).
\newblock Positive definite metric spaces.
\newblock Preprint, \arXiv{1012.5863}.

\bibitem[MT]{MT} J. Melleray and T. Tsankov (2011).
Generic representations of
abelian groups and extreme amenability. Preprint, \arXiv{1107.1698v1}.

\bibitem[Mo]{Mo} J.T. Moore (2011).
Amenability and Ramsey theory. Preprint, 
\arXiv{1106.3127v4}.

\bibitem[NR]{NR}
J. Ne{\v{s}}et{\v{r}}il and{} V. R{\"o}dl (1978).
\newblock On a probabilistic graph-theoretical method.
\newblock {\em Proc. Amer. Math. Soc.} {\bf 72}, 417--421.

\bibitem[N1]{N1} L. Nguyen Van Th\'e (2010).
Structural Ramsey theory of metric
spaces and topological dynamics of isometry groups, {\it Memoirs of
the Amer. Math. Soc.} {\bf 206}, No. {\bf 968}.

\bibitem[N2]{N2} L. Nguyen Van Th\'e (2012).
More on the Kechris-Pestov-Todorcevic correspondence: Precompact
expansions. Preprint,
\arXiv{1201.12708v1}.

\bibitem[Pe1]{Pe1} V. Pestov (1998). On free actions, minimal flows and a
problem by Ellis, {\it Trans. Amer. Math. Soc.} {\bf 350} (10), 4149--4165.

\bibitem[Pe2]{Pe2} V. Pestov (2002). Ramsey-Milman phenomenon, Urysohn
metric spaces, and extremely amenable groups, {\it Israel J. Math.} {\bf 127},
317--357.

\bibitem[Sch]{Sch} 
I.J. Schoenberg (1938).
\newblock Metric spaces and positive definite functions.
\newblock {\em Trans. Amer. Math. Soc.} {\bf 44}, 522--536.

\bibitem[So]{So} S. Solecki (2005). Extending partial isometries, {\it
Israel J. Math.} {\bf 150}, 315--332.

\bibitem[W]{W} B. Weiss (2012). Minimal models for free actions, {\it
Contemp. Math.} {\bf 567}, 249--264.

\end{thebibliography}
\end{document}